\documentclass[10pt,reqno]{amsart}

\usepackage[utf8]{inputenc}
\usepackage[T1]{fontenc}

\usepackage{amsmath, amssymb, amsthm}
\usepackage{geometry}
\usepackage{enumitem}
\usepackage{hyperref}
\usepackage{marginnote}
\usepackage{xcolor}
\usepackage{upgreek}
\usepackage{tikz}
\usetikzlibrary{calc}
\usepgflibrary {shadings} 
\usetikzlibrary{intersections}
\usetikzlibrary{patterns}
\usepackage{mathtools}
\usepackage{verbatim}

\numberwithin{equation}{section}

\theoremstyle{plain}
\newtheorem{theorem}{Theorem}
\newtheorem{lemma}{Lemma}[section]
\newtheorem{proposition}{Proposition}[section]

\theoremstyle{definition}

\theoremstyle{remark}

\newcommand{\N}{\mathcal{N}}

\newcommand{\C}{\mathbb{C}}
\newcommand{\R}{\mathbb{R}}

\newcommand{\round}[1]{\left({#1}\right)}
\newcommand{\squared}[1]{\left[{#1}\right]}
\newcommand{\braces}[1]{\left\lbrace{#1}\right\rbrace}
\newcommand{\abs}[1]{\left\lvert{#1}\right\rvert}
\newcommand{\norm}[1]{\left\lVert{#1}\right\rVert}
\newcommand{\bracket}[1]{\left\langle{#1}\right\rangle}

\newcommand{\RE}[1]{\mathrm{Re}\left({#1}\right)}

\newcommand{\sdf}[1]{\mathsf{#1}}
\newcommand{\twovec}[2]{\begin{pmatrix}{#1}\\{#2}\end{pmatrix}}

\newcommand{\Pssi}[1]{\boldsymbol{\Psi}}

\newcommand{\eps}{\varepsilon}
\newcommand{\J}{\mathcal{J}}

\newcommand{\vphi}{\varPhi}
\newcommand{\U}{\mathcal{U}}
\newcommand{\V}{\mathcal{V}}
\newcommand{\Ric}{\mathrm{Ric}}

\begin{document}

\title[Solutions for magnetic Ginzburg--Landau equations]{Solutions of the Ginzburg--Landau equations concentrating on codimension-2 minimal submanifolds}

    \author[M. Badran]{Marco Badran}
    
	\address{Department of Mathematical Sciences, University of Bath, Bath, BA2 7AY, United Kingdom. }
 
    \email{mb2747@bath.ac.uk}	

	\author[M. del Pino]{Manuel del Pino}
	
	\address{Department of Mathematical Sciences, University of Bath, Bath, BA2 7AY, United Kingdom.}
	
	\email{mdp59@bath.ac.uk}

\maketitle

\begin{abstract}
We consider the magnetic Ginzburg-Landau equations in a compact manifold $N$
\begin{equation*}
\begin{cases}
    -\eps^2\Delta^Au =
  \frac{1}{2}(1-|u|^{2})u,\\
    \eps^2d^{*}dA=\bracket{\nabla^Au,iu},
    \end{cases}
\end{equation*}
formally corresponding to the Euler-Lagrange equations for the energy functional
\begin{equation*}
     E(u,A)=\frac{1}{2}\int_{N}\eps^2|\nabla^Au|^{2}+\eps^4|dA|^{2}+\frac{1}{4}(1-|u|^{2})^{2}.
\end{equation*}
Here $u:N\to \C$ and $A$ is a 1-form on $N$.
Given a codimension-2 minimal submanifold  $M\subset N$ which is also oriented and non-degenerate, we construct a solution  $(u_\eps,A_\eps)$ such that $u_\eps$ has a
zero set consisting of a smooth surface close to $M$. Away from $M$ we have
\begin{equation*}
	u_\eps(x)\to\frac{z}{|z|},\quad A_\eps(x)\to\frac{1}{|z|^2}(-z_2dz^1+z_1dz^2),\quad x=\exp_y(z^\beta\nu_\beta(y)).
\end{equation*}
as $\eps\to 0$, for all sufficiently small $z\ne 0$ and $y\in M$. Here, $\{\nu_1,\nu_2\}$ is a normal frame for $M$ in $N$. This improves a recent result by De Philippis and Pigati \cite{DePhilippis-Pigati2022} who built a solution for which the concentration phenomenon holds in an  energy, measure-theoretical sense.  
\end{abstract}

\section{Introduction}
Let $n\geq 4$. We consider for $\eps>0$ the magnetic Ginzburg-Landau energy on a smooth, closed $n$-dimensional manifold $N$,
\begin{equation}\label{energy}
	E_\eps(u,A)=\frac12\int_N\eps^2|\nabla^Au|^2+\eps^4|dA|^2+\frac14\round{1-|u|^2}^2
\end{equation}
where $u:N\to\C$ is a complex valued function and $A\in \Omega^1(N)$ is the magnetic potential. Here $d$ denotes the exterior differential and $\nabla^A=d-iA$ is the covariant derivative. Explicitly, if $\{x^1,\dots,x^n\}$ are local coordinates on $N$ and $g$ is the metric, then we write $A=A_jdx^j$ and 
\begin{align*}
	|\nabla^Au|^2&=\sum_{j,k=1}^ng^{ij}(\partial_ju-iA_ju)(\partial_k\bar u+iA_k\bar u),\\
	|dA|^2&=\frac12\sum_{j,k,s,t=1}^ng^{ks}g^{jt}(\partial_kA_j-\partial_jA_k)(\partial_sA_t-\partial_tA_s).
\end{align*}
The corresponding Euler-Lagrange system is given by
\begin{equation}\label{equations}
		\begin{cases}
			-\eps^2\Delta^Au=\frac12(1-|u|^2)u,\\
			\eps^2d^*dA=\langle\nabla^Au,iu\rangle
		\end{cases}\quad\text{on }N,
\end{equation}
which can be also written as $S_\eps(u,A)=0$, where $S_\eps=E'_\eps$. We denoted with $-\Delta^A=(\nabla^A)^*\nabla^A$ the connection laplacian. The operators in \eqref{equations} explicitly read
\begin{align*}
    -\Delta^Au&=-\frac{1}{\sqrt{\det g}}(\partial_j-iA_j)\left[\sqrt{\det g}g^{jk}(\partial_k-iA_k)u\right],\\
    d^*dA&=-\frac{1}{\sqrt{\det g}}g_{jk}\partial_{i}\left(\sqrt{\det g}g^{li}g^{tk}(\partial_l A_{t}-\partial_tA_l)\right)dx^{j}.
\end{align*}

Energy \eqref{energy} appears in the classical theory of superconductivity by Ginzburg and Landau when $N=\R^n$, $n=2,3$. In \cite{Taubes1980_1} were found solutions of \eqref{equations} in the planar case $\R^2$ with isolated zeros (vortices) of $u$, and in particular in \cite{Berger-Chen1989} was found a solution of \eqref{equations} with $\eps=1$ in $\R^2$ with a degree 1 radial symmetry, namely $U_0=(u_0,A_0)$ where
\begin{equation}\label{building-block}
	u_0(\zeta)=f(r)e^{i\theta},\quad A_0(\zeta)=a(r)d\theta,\quad \zeta=re^{i\theta}.
\end{equation} 
The functions $f(r)$ and $a(r)$ are positive solutions to the system of ODE 
\begin{equation}
\begin{cases}
		-f''-\frac{f'}{r}+\frac{(1-a)^2f}{r^2}-\frac{1}{2}f(1-f^2)=0,\\
-a''+\frac{a'}{r}-f^2(1-a)=0
\end{cases}
\quad\text{in }(0,+\infty)
\end{equation}
with $f(0)=a(0)=0$.
This is the unique solution of \eqref{equations} with $\eps=1$ and exactly one zero with degree 1 at the origin. Also, $U_0$ is linearly stable as established in \cite{Gustafson-Sigal2000, Stuart1994}. In addition
\begin{equation*}
	f(r)-1=O(e^{-{r}}),\quad a(r)-1= O(e^{-r}),\quad\text{as }r\to \infty,
\end{equation*}
see \cite{Berger-Chen1989, Plohr1981}. 

In this paper we look for of solutions $U_\eps=(u_\eps,A_\eps)$ of \eqref{equations} concentrating in the limit $\eps\to 0$ around a given codimension 2 minimal submanifold $M\subset N$, in the form of $\eps$-scalings of \eqref{building-block}, with $z$ representing a normal coordinate of $M$ inside $N$. More precisely, let $\{\nu_1,\nu_2\}$ be an orthonormal basis for $T^\perp M$. We describe a tubular neighbourhood of $M$ in $N$ by local Fermi coordinates 
\begin{equation}\label{tub neigh}
	x=X(y,z)=\exp_y(z^\beta\nu_\beta(y)),\quad y\in M,\ \  |z|<\tau
\end{equation}
for some $\tau>0$. Then, the solution $U_\eps(x)=(u_\eps(x),A_\eps(x))^T$ behaves as 
\begin{equation*}
	u_\eps(x)\approx f\round{\frac{z}\eps}\frac{z}{|z|},\quad A_\eps(x)\approx a\round{\frac{z}\eps}\frac{1}{|z|^2}(-z_2dz^1+z_1dz^2).
\end{equation*}
In particular it satisfies
\begin{equation*}
	u_\eps(x)\to\frac{z}{|z|},\quad A_\eps(x)\to\frac{1}{|z|^2}(-z_2dz^1+z_1dz^2),\quad \eps\to 0.
\end{equation*}
We consider a closed, $n$-dimensional manifold $N^n$ and a closed $n-2$-dimensional minimal submanifold $M^{n-2}\subset N^n$. We  say  that a minimal manifold $M\subset N$ is \emph{admissible} if
\begin{itemize}
	\item[(H)] - $M$ is the
	boundary of a $n-1$-dimensional, oriented, 
	embedded submanifold $B^{n-1}\subset N^n$.  
\end{itemize} 
Condition (H) is adapted from \cite{Colinet-Jerrard-Sternberg2021}. We also require that $M$ is non-degenerate, in the sense that the Jacobi operator has trivial bounded kernel, namely
\begin{equation}\label{nondeg}
    h\in L^\infty(M),\ \J[h] = 0\implies h=0.
\end{equation}
We recall that the Jacobi operator is the second variation of the area functional around $M$, namely $\J[h]=(\J^1[h],\J^2[h])^T$ where
\begin{equation*}
	\J^\gamma[h]=\Delta_Mh^\gamma+k_{\ell\beta}k_{\ell\gamma}h^\beta+\Ric(\nu_\gamma,\nu_\beta)h^\beta,\quad \gamma=1,2.
\end{equation*}
Now we state the main result. 
\begin{theorem}\label{main theorem}
	Let $N$ be a closed $n$-dimensional manifold and let $M\subset N$ be an admissible, non-degenrate, codimension-2 minimal submanifold. Then there is $\delta>0$ such that for $\sigma\in(0,1)$ and all sufficiently small $\eps>0$ there exists a solution $(u_\eps,A_\eps)$ to \eqref{equations} which as $\eps\to 0$ satisfies 
	\begin{equation}\label{local info}
	\begin{split}
		u_\eps(x)&=u_0\round{\frac{z-\eps^2h_0(y)}{\eps}}+O\round{\eps^2 e^{-\frac{\sigma|z|}{\eps}}},\\
		A_\eps(x)&=A_0\round{{\frac{z-\eps^2h_0(y)}{\eps}}}+O\round{\eps^2 e^{-\frac{\sigma|z|}{\eps}}},
	\end{split}\qquad |z|<\delta,
	\end{equation}
	for all points $x=X(y,z)$ of the form \eqref{tub neigh} and where $h_0$ is a smooth function on $M$, explicitly characterized in \eqref{eq h0} below.
	Moreover, $|u_\eps|\to 1$ uniformly on compact subsets of $N\setminus M$. The zero set of $u_\eps$
consists of a smooth codimension-2 manifold,  $\eps^2$-close to $M$.
\end{theorem} 

In the interesting paper \cite{DePhilippis-Pigati2022}, De Philippis and Pigati consider a setting similar to that of Theorem \ref{main theorem} and construct a solution exhibiting energy concentration near $M$. Their proof does not yield the refined local description \eqref{local info} nor the manifold character of the zero set.

Their approach, variational in nature, is different from ours. Our strategy relies on the construction of of an accurate approximation to a solution, the use of linearization and the so-called outer-inner gluing method to set up a suitable fixed point formulation. We have recently used this approach in a similar context in $\R^4$ 
obtaining concentration for a special class of non-compact minimal surfaces embedded in $\R^3$ which includes a catenoid and the Costa-Hoffman-Meeks minimal surfaces \cite{Badran-delPino2022}. 
Liu, Ma, Wei and Wu \cite{Liu-Ma-Wei-Wu2021} have obtained the existence of a solution in $\R^4$ with precise asymptotics and concentration on a special non-compact codimension-2  symmetric minimal  manifolds discovered by Arezzo and Pacard \cite{Arezzo-Pacard2003}. See also \cite{Davila-delPino-Medina-Rodiac2022} for a recent construction in $\R^3$ of interacting helicoidal vortex filaments in Ginzburg-Landau without magnetic field. 

\medskip 
In the related Allen-Cahn equation, analogues of these  constructions via gluing methods appear in 
\cite{Pacard-Ritore2003, delPino-Kowalczyk-Wei2011,delPino-Kowalczyk-Wei2013, delPino-Pacard-Wei2015}. 
In the context of the Ginzburg-Landau equation, with or without magnetic fields related variational constructions have been performed in \cite{Colinet-Jerrard-Sternberg2021, Montero-Sternberg-Ziemer2004, Riviere1996, Jerrard-Sternberg2009, Lin-Riviere1999,Bethuel-Brezis-Orlandi2001, Contreras-Jerrard2017}

\medskip
Theorem \ref{main theorem} relates to a program of construction of codimension-2 minimal submanifolds as limits of critical points of scaling energies  Ginzburg-Landau type, see Pigati-Stern \cite{Pigati-Stern2021}.   This program has already been successful in the codimension-1 case, see \cite{Guaraco2018, Gaspar-Guaraco2018,  Bellettini-Wickramasekera2020,Chodosh-Mantoulidis2020} where the Allen-Cahn energy is used as an  alternative to the Almgren-Pitts min-max theory to construct minimal surfaces.

\medskip
We now recall some facts that will be relevant in the rest of the paper.
As already mentioned, we work with pairs $W=(u,A)^T$ given by a complex-valued function  $u$ and a 1-form $A$ defined on $N$. We define an $\eps$-dependent $L^2$-inner product on these pairs by setting 
\begin{align*}
	\langle W_1,W_2\rangle_\eps&=\int_N\twovec{u_1}{A_1}{\cdot_\eps}\twovec{u_2}{A_2}\\
	&=\int_N\langle u_1,u_2\rangle+\eps^2A_1\cdot A_2
\end{align*}
where $\langle u_1,u_2\rangle=\RE{u_1\bar{u}_2}$ and $A_1\cdot A_2=g^{ij}(A_1)_i(A_2)_j$, being $g$ the metric on $N$. Recall that energy \eqref{energy} is invariant under $U(1)$-gauge transformations $\mathrm{G}_\gamma$, given by 
\begin{equation*}
	\mathrm{G}_\gamma\twovec{u}{A}=\twovec{ue^{i\gamma}}{d\gamma},
\end{equation*}
for any smooth function $\gamma$ on $N$. This invariance creates an infinite dimensional part of the kernel of the linearized operator $S'_\eps(W)$, which we will denote with $Z_{W,g}$. Such part of the kernel is generated by the gauge-zero modes, namely $Z_{W,g}=\mathrm{span}_\gamma\{\Theta_W[\gamma]\}$, with 
\begin{equation*}
	\Theta_W[\gamma]=\twovec{iu\gamma}{d\gamma}
\end{equation*}
and where again $\gamma$ varies among all smooth functions defined on $N$.
Observe that we can characterize the orthogonality to such kernel with respect to the inner product $\langle\cdot,\cdot\rangle_\eps$ by using the adjoint of $\Theta_W$,
\begin{equation*}
	\Phi\in Z_{W,g}^\perp\iff \Theta_W^*[\Phi]=0,
\end{equation*}
where, if $\Phi=(\phi, \omega)^T$
\begin{equation*}
	\Theta_W^*[\Phi]=\eps^2d^*\omega+\langle\phi,iu\rangle.
\end{equation*}
With this notation, the linearized operator can be expressed as 
\begin{equation}\label{decomposition}
	S_\eps'(W)=L^\eps_W-\Theta_W\Theta_W^*
\end{equation}
where 
\begin{equation}\label{LW}
	L^\eps_W[\Phi]=\twovec{-\eps^2\Delta^A\phi-\frac12(1-3|u|^2)\phi+2i\eps^2\nabla^Au\cdot\omega}{-\eps^2\Delta\omega+|u|^2\omega+2\langle\nabla^Au,i\phi\rangle}
\end{equation}
is an elliptic operator, being $-\Delta\omega=(d^*d+dd^*)\omega$ the Hodge laplacian on 1-forms. We remark that the natural space in which $L_W$ is defined and smooth is the space $H_W^1(N)$ of pairs $\Phi=(\phi,\omega)^T$ such that 
\begin{equation}\label{H1N}
	\|\Phi\|_{H^1_W(N)}\coloneqq\|\nabla^A\phi\|_{L^2(N)}+\|\phi\|_{L^2(N)}+\|\nabla\omega\|_{L^2(N)}+\|\omega\|_{L^2(N)}<\infty
\end{equation}
where $\nabla\omega$ is the Levi-Civita connection applied to the 1-form $\omega$. 
 Decomposition \eqref{decomposition} is meaningful because 
\begin{equation}\label{properties ThetaW}
	\mathrm{range}(\Theta_W\Theta_W^*)\subset Z_{W,g},\quad \ker(\Theta_W\Theta_W^*) = Z_{W,g}^\perp.
\end{equation}
Finally we define, for $W=(u,A)^T$ the gradient-like operator 
\begin{equation}\label{gradient-like}
	\nabla_W\twovec{\phi}{\omega}=\twovec{\nabla^A\phi}{d\omega+d^*\omega}
\end{equation}
and the corresponding Laplacian-like operator $-\Delta_W=\nabla^*_W\nabla_W$, which explicitly reads
\begin{equation*}
	-\Delta_W\twovec{\phi}{\omega}=\twovec{-\Delta^A\phi}{-\Delta\omega}.
\end{equation*}
This notation allows us to rewrite \eqref{LW} as 
\begin{equation}\label{LW compact}
	L^\eps_W[\Phi]=-\eps^2\Delta_W\Phi+\Phi+T^\eps_W\Phi
\end{equation}
where
\begin{equation*}
	T^\eps_W\twovec{\phi}{\omega}=\twovec{-\frac12(1-|u|^2)\phi+2i\eps^2\nabla^Au\cdot\omega}{-(1-|u|^2)\omega+2\langle\nabla^Au,i\phi\rangle}.
\end{equation*}

\medskip
We will use extensively the case in which the space is the whole plane $\R^2$. There is a rich literature devoted to this problem, see for instance \cite{Taubes1980_1, Taubes1980_2}. 
Recall the degree 1 solution $U_0$ defined in \eqref{building-block}; it is convenient to denote 
\begin{equation}\label{lstr}
	\sdf{L}\coloneqq L_{U_0}=S'(U_0)-\Theta_{U_0}\Theta_{U_0}^*
\end{equation}
where $L_{U_0}\coloneqq L_{U_0}^1$.
We recall that
\begin{equation*}
	Z_{U_0}=\ker S'({U_0})=Z_{{U_0},t}\oplus Z_{{U_0},g}
\end{equation*}
where $Z_{U_0,t}=\mathrm{span}\{\sdf{V}_1,\sdf{V}_2\}$, being
\begin{equation}\label{expression sdfVj}
	\sdf{V}_1=\twovec{f'}{\frac{a'}{r}dt^2},\quad \sdf{V}_2=\twovec{if'}{-\frac{a'}{r}dt^1}.
\end{equation}
Also, as shown by Stuart \cite{Stuart1994}, it holds the coercivity estimate
\begin{equation}\label{coercivity}
	\langle \sdf{L}[\Phi],\Phi\rangle_{L^2}\geq c\norm{\Phi}^2_{H^{1}_{U_0}},\quad\forall\Phi\in Z_{{U_0},t}^\perp
\end{equation}
for some $c>0$, where 
\begin{equation*}
	\norm{\Phi}^2_{H^{1}_{U_0}}=\norm{\nabla^A\phi}^2_{L^2}+\norm{\phi}^2_{L^2}+\norm{\nabla\omega}^2_{L^2}+\norm{\omega}^2_{L^2}
\end{equation*}
is the covariant norm. This coercivity estimate along with Lax-Milgram theorem imply the validity of the following existence result.
\begin{lemma}\label{invertibility lstr}
	Let $\Psi\in L^2(\R^2)$ satisfy
	\begin{equation}\label{lemma1 ort cond}
		\int_{\R^2}\Psi\cdot\sdf{V}_\alpha dt=0,\quad \alpha=1,2,
	\end{equation} 
	where $\sdf{V}_\alpha$ are as in \eqref{expression sdfVj}. Then there exists a unique solution $\Phi\in H^1_{U_0}(\R^2)$ to 
	\begin{equation*}
		\sdf{L}[\Phi]=\Psi
	\end{equation*}
	 satisfying 
	 \begin{equation*}
	 	\int_{\R^2}\Phi\cdot\sdf{V}_\alpha dt=0,\quad \alpha=1,2,
	 \end{equation*}
	 and
	\begin{equation*}
		\|\Phi\|_{H^1_{U_0}(\R^2)}\leq C\|\Psi\|_{L^2(\R^2)}
	\end{equation*}
	for some $C>0$.
\end{lemma}

\subsection{The role of gauge invariance in the linearized equation}\label{role gauge}

Here we recall briefly the content of \S2.2 of \cite{Badran-delPino2022}. 
To prove Theorem \ref{main theorem} we first construct an approximate solution $W$, whose profile close to the manifold resembles the behaviour of $U_0$, then we look for true solutions as perturbations $\vphi$ of the first approximation $W$, namely we solve 
\begin{equation}\label{perturbed eq Seps}
	S_\eps(W+\vphi)=0.
\end{equation}
By denoting $N_\eps(\vphi)=S_\eps(W+\vphi)-S_\eps(W)-S'_\eps(W)[\vphi]$, we formulate equation \eqref{perturbed eq Seps} as
\begin{equation}\label{perturbed eq Seps expanded}
	S'_\eps(W)[\Phi]=-S_\eps(W)-N_\eps(\Phi).
\end{equation}
Our task is to develop an invertibility theory for $S'_\eps(W)$ in order to rephrase \eqref{perturbed eq Seps expanded} as a fixed point problem and then solve it using the small Lipschitz character of $N_\eps$. The infinite-dimensionality of the kernel of $S'_\eps(W)$ could be an obstruction to a convenient invertibility theory. We address this issue by projecting \eqref{perturbed eq Seps} onto both the gauge kernel and its orthogonal
\begin{align}
	\pi_{Z_{W,g}^\perp}S_\eps(W+\vphi)&=0,\label{projection}\\
	\pi_{Z_{W,g}}S_\eps(W+\vphi)&=0.\label{projection-orthogonal}
\end{align}
Recalling the decomposition \eqref{decomposition} and the fact that, by self-adjointness of $S_\eps'(W)$, 
\begin{equation}\label{selfmap1}
	L_W[Z_{W,g}^\perp]=S_\eps'(W)[Z_{W,g}^\perp]\subset Z_{W,g}^\perp
\end{equation}
we infer that by choosing $\vphi\in Z_{W,g}^\perp$ equation \eqref{projection} reads 
\begin{equation}\label{projected eq extended}
	L_W[\vphi]=-S_\eps(W)-\pi_{Z_{W,g}^\perp}N(\vphi).
\end{equation}
In this way, we reduce the task of inverting $S'_\eps(W)$ to that of inverting $L_W$, which doesn't have the infinite gauge-kernel issue. Moreover, from \eqref{selfmap1} and the fact that 
\begin{equation*}
	L_W[Z_{W,g}]=\Theta_W\Theta_W^*[Z_{W,g}]\subset Z_{W,g}.
\end{equation*}
we get that, if an inverse $L_W^{-1}$ exists, then 
\begin{equation}\label{self-map}
	L_W^{-1}[Z_{W,g}^\perp]\subset Z_{W,g}^\perp.
\end{equation}
This fact implies that formulating \eqref{projected eq extended} as a fixed point by applying $L_W^{-1}$ the right-hand side will be well defined as a map $ Z_{W,g}^\perp\mapsto  Z_{W,g}^\perp$. We will manage to apply this principle and find a solution $\vphi$ belonging to a sufficiently small ball in a suitable norm. Now we use the fact that, as pointed out in \cite{Ting-Wei2013}, if \eqref{projection-orthogonal} holds for a $\vphi$ that is sufficiently small with respect to $W$ then also \eqref{projection} holds. To see this, suppose the opposite, namely the existence of a function $\gamma$ such that 
\begin{equation*}
	\pi_{Z_{W,g}^\perp}S_\eps(W+\vphi)=\Theta_W[\gamma].
\end{equation*}
Let $W=(u,A)$ and $\vphi=(\varphi,\omega)$, then we have 
\begin{align*}
	0&=\int_NS_\eps(W+\vphi)\cdot_\eps\Theta_{W+\vphi}[\gamma]\\
	&=\int_N\Theta_W[\gamma]\cdot_\eps\Theta_{W+\vphi}[\gamma]\\
	&=\int_N\gamma\cdot\squared{(-\eps^2\Delta+|u|^2+\langle u,\varphi\rangle)\gamma}
\end{align*}
and since $-\eps^2\Delta+|u|^2+\langle u,\varphi\rangle$ is a positive operator given that $\varphi$ is small with respect to $u$, it must be $\gamma=0$.

\section{The first approximation}

\begin{figure}
	\centering
	\begin{tikzpicture}
	\coordinate (O) at (1,2.5);
		\shadedraw[thick, top color=blue!70, bottom color=blue!30, fill opacity=.3] (2.5,3.75) to[in=175, out=45] (5,4) to[out=-5,in=160] (9,4.1) to[out=-20, in=20]  (9,4) to[out=45, in=-25] (8,4.5) to[out=155,in=-10] (6,5.4)  to[out=170,in=55] (2.5,3.75);
		
		\draw[black, thin, fill=gray, fill opacity=.2] (9,4) to[out=45, in=-25] (8,4.5) to[out=155,in=-10] (6,5.4)  to[out=170,in=55] (2.5,3.75) to[out=240, in=185] (4.1,3.3) to[out=5, in=200] (9,4);
		
		\shadedraw[thick, top color=blue, bottom color=blue!30, fill opacity=.3] (2.5,3.75) to[in=175, out=45] (5,4) to[out=-5,in=160] (9,4.1) to[out=-20, in=20]  (9,4) to[out=200, in=5] (4.1,3.3) to[out=185, in=240] (2.5,3.75);
		
		\draw[orange, ultra thick, draw opacity =1] (2.5,3.75) to[out=240, in=185] (4.1,3.3) to[out=5, in=200] (9,4);
		
		\draw[orange, ultra thick] (2.5,3.75) to[in=175, out=45] (5,4) to[out=-5,in=160] (9,4.1) to[out=-20, in=20]  (9,4);
		
		\draw [red,->, thick]  (3,4.07) -- ++(50:.6cm);
		\draw [blue,->, thick]  (3,4.07) -- ++(140:.6cm);
		
		\draw[->,thick] (O) -- ++ (0:2.5cm);
		\draw[->,thick] (O) -- ++ (90:3.5cm);
		\draw[->,thick] (O) -- ++ (215:1.5cm);
		
		\node[right] at (5,5) {$B$};
		\node[right] at (3.4,4.5) {$\nu_2$};
		\node[left] at (2.6,4.5) {$\nu_1$};
		\node[orange] at (8.5,3.5) {$M$};
		\node at (1.5,6) {$N$};
	\end{tikzpicture}
	\caption{A representation of $M$ as the boundary of an oriented manifold $B$ in $N$. Assumption (H) determines the normal fields $\{\nu_1,\nu_2\}$. }
	\label{fig admissibility}
\end{figure}
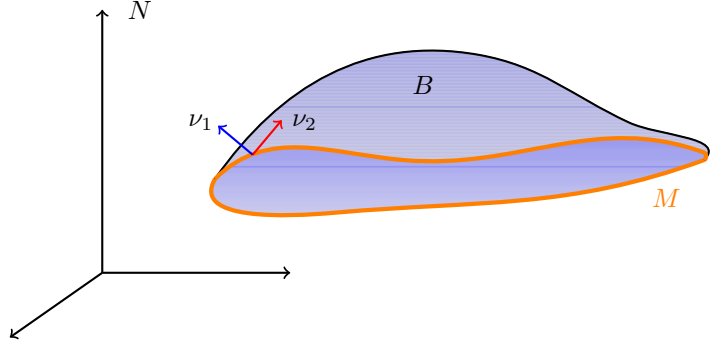

The setting of Theorem \ref{main theorem} allows us to choose in a standard way a basis for the normal bundle $T^\perp M$ for the immersion $M\subset N$. 
Recall that, by the admissibility hypothesis (H), $M=\partial B$ where $B$ is oriented.
Let $\nu_2$ be the normal field of the immersion $B\subset N$ (which exists because of the orientability assumption) and let $\nu_1$ the vector field of $TB$ restricted to $M$ which is normal to $TM$ (we can assume that $\nu_1$ points inwards $B$). In this way $\{\nu_1,\nu_2\}$ is a basis of $T^\perp M$, see Figure \ref{fig admissibility}. In particular, given a basis $\{e_1,\dots,e_{n-2}\}$ for $TM$ then $\{e_1,\dots,e_{n-2},\nu_1,\nu_2\}$ is a basis of $TN$ defined on $M$. 

In what follows letters $i,j,k,\ldots$ are used for tangential coordinates to the manifold $M$, while $\alpha,\beta,\gamma,\ldots$ denote coordinates in the normal directions. We use $a,b,c\ldots$ 
 to indicate all coordinates at once. More precisely,
 \begin{equation*}
 	1\leq i,j,k,\ldots\leq n-2,\quad 1\leq \alpha,\beta,\gamma,\ldots\leq 2,\quad 1\leq a,b,c,\ldots\leq n.
 \end{equation*}
We describe a tubular neighbourhood of $M$ in $N$ with coordinates defined by the exponential map,
\begin{equation}\label{standard Fermi}
	x=X(y,z) = \exp_y(z^\beta\nu_\beta(y)),\quad (y,z)\in M\times B(0,\tau)
\end{equation}
where $B(0,\tau) \subset \R^2$ and $\tau$ is sufficiently small. 
We consider the change of coordinates
\begin{equation}\label{change of coordinates}
	z^\beta =\eps(t^\beta+h^\beta(y))
\end{equation}
for some pair $h=(h^1,h^2)$ of functions defined on $M$. This produces new coordinates 
\begin{equation}\label{Xh}
	x=X_h(y,t)=\exp_y\round{\eps(t^\beta+h^\beta(y))\nu_\beta(y)}
\end{equation}
defined on
\begin{equation*}
	\mathcal{O}_h=\braces{(y,t)\ :\ y\in M,\ |t+h(y)|<\tau/\eps}
\end{equation*}
and we set $\N_h=X_h(\mathcal{O}_h)$. 
\begin{figure*}
	\begin{tikzpicture}

\draw [black, thick] plot [smooth] coordinates {(0,4) (2,4) (3,5.5) (5,5.5) (5.5,5)  (4,3.5) (5,3) (7,3) (8.5,3.9) (9,4)};

\coordinate (U) at (-1.5,3.5);
\node[right] at (9,4) {$M$};
\draw [thick,blue,->] (U) -- ++ (90:1.1cm);
\draw [thick,red,->] (U) -- ++ (0:1.1cm);
\node[red, above right] at (-0.9,3.5) {$\nu_1$};
\node[blue, below left] at (-0.9,4.7) {$\nu_2$};

\coordinate (A) at (1,2.66);
\def\Al{0.95cm}
\def\Ah{2.1cm}
\def\Aangle{30}
\def\AVangle{90}
\coordinate (Acent) at (1.43,3.91);

\draw[fill=white, fill opacity=0.9,draw opacity=0] plot coordinates {(A) (Acent) (1, 4.5)};
\draw[fill=gray, fill opacity=0.2, draw opacity=0.2] (A) -- ++(\Aangle:\Al) -- ++(\AVangle:\Ah) -- ++ (\Aangle+180:\Al) -- ++ (-\AVangle:\Ah);
\filldraw[red] (Acent) circle (0.03);

\draw [blue,->] (Acent) -- ++ (\AVangle:0.5cm);
\draw [red,->] (Acent) -- ++ (\Aangle+180:0.4cm);

\coordinate (B) at (4.8,3.64);
\def\Bl{1.1cm}
\def\Bh{1.8cm}
\def\Bangle{10}
\def\BVangle{110}
\coordinate (Bcent) at (5.05,4.54);

\draw[fill=white, fill opacity=0.9,draw opacity=0] plot coordinates {(5.63,4.5) (Bcent) (5, 5.47) (5.26, 5.52)};
\draw[fill=gray, fill opacity=0.2, draw opacity=0.2] (B) -- ++(\Bangle:\Bl) -- ++(\BVangle:\Bh) -- ++ (\Bangle+180:\Bl) -- ++ (\BVangle+180:\Bh);
\filldraw[red] (Bcent) circle (0.03);

\draw [blue,->] (Bcent) -- ++ (\BVangle:0.5cm);
\draw [red,->] (Bcent) -- ++ (\Bangle+180:0.4cm);

\coordinate (C) at (7,2.5);
\def\Cl{1.1cm}
\def\Ch{1.8cm}
\def\Cangle{-10}
\def\CVangle{80}
\coordinate (Ccent) at (7.69,3.34);

\draw[fill=white, fill opacity=0.9,draw opacity=0] plot coordinates { (8.2,3) (Ccent) (7.4,3.34) (8.39, 4.08) };
\draw[fill=gray, fill opacity=0.2, draw opacity=0.2] (C) -- ++(\Cangle:\Cl) -- ++(\CVangle:\Ch) -- ++ (\Cangle+180:\Cl) -- ++ (\CVangle+180:\Ch);
\filldraw[red] (Ccent) circle (0.03);

\draw [blue,->] (Ccent) -- ++ (\CVangle:0.5cm);
\draw [red,->] (Ccent) -- ++ (\Cangle:0.4cm);

	\end{tikzpicture}
	\caption{The normal frame around $M$.}
\end{figure*}
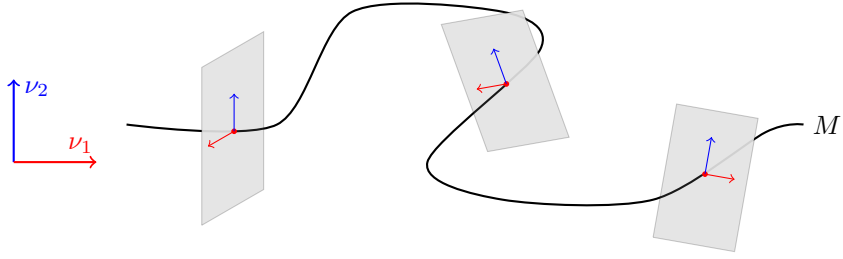
With this choice of coordinates, we define the first local approximation $W_0$ to a solution of Problem \eqref{equations} as 
\begin{equation}\label{new coordinates}
	W_0(x)=U_0(t),\quad x=X_h(y,t)
\end{equation}
with $U_0(t)=(u_0,A_0)^T$ defined in \eqref{building-block}. 
Below, we shall estimate the error of approximation $S_\eps(W_0)$, where 
\begin{equation}\label{S_eps}
	S_\eps(u,A)=\twovec{-\eps^2\Delta^Au-\frac12(1-|u|^2)u}{\eps^2d^*dA-\langle\nabla^Au,iu\rangle}.
\end{equation}
It is convenient to express the operators in $S_\eps$ in coordinates $(y,t)$ in \eqref{Xh}, in particular the operators $-\Delta_N^A$ and $d^*_Nd_N$. First, we express them in coordinates $(y,z)$ in \eqref{standard Fermi}. Standard computations, see for instance \cite[Appendix A.1]{Badran-delPino2022}, show that the metric in $N$, can be expressed as 
\begin{equation}\label{gN}
	g_N(y,z)=\begin{pmatrix}g_{M_z}(y) & 0\\
0 & I
\end{pmatrix}
\quad\text{on }M\times B(0,\tau)
\end{equation}
where $g_{M_z}(y)$ is the metric of $N$ restricted to 
\begin{equation*}
	M_z=\braces{\exp_y(z^\beta\nu_\beta(y))\ :\ y\in M}.
\end{equation*}
Before proceeding, we make explicit our assumption on $h$. First, we set in a standard way, for $0<\gamma <1$, 
\begin{equation*}
	\|\varphi\|_{C^{0,\gamma}(\Lambda)}=\|\varphi\|_{L^\infty(\Lambda)}+\sup_{p\in \Lambda}[\varphi]_{\gamma,B(p,1)}
\end{equation*}
where $\varphi$ is a function defined on a manifold $\Lambda$, $B(p,1)$ is the geodesic ball of radius 1 around $p\in \Lambda$, and
\begin{equation*}
	[\varphi]_{\gamma,X}=\sup_{\substack{x,y\in X\\x\ne y}}\frac{|\varphi(x)-\varphi(y)|}{|x-y|^\gamma}.
\end{equation*}
We also set, for $k\geq 1$, 
\begin{equation*}
	\|\varphi\|_{C^{k,\gamma}(\Lambda)} = \|\varphi\|_{C^{0,\gamma}(\Lambda)}+\|D^k\varphi\|_{C^{0,\gamma}(\Lambda)}.
\end{equation*}
With this notation, we state our hypothesis on $h$: we assume that for some fixed $K>0$ 
\begin{equation}\label{assumption h}
	\|h\|_{C^{2,\gamma}(M)}\leq K\eps.
\end{equation}
 
\subsection{The operators in coordinates}
The following expressions for the differential operators hold. Given any 1-form $\omega=\omega_adx^a\in \Omega^1(N)$ the covariant Laplacian $-\Delta_N^\omega$ is given by 
\begin{equation*}
	-\Delta_N^\omega\phi=-\frac{1}{\sqrt{\det g_N}}\partial^\omega_a\round{\sqrt{\det g_N}g_N^{ab}\partial^\omega_b\phi}
\end{equation*}
where $\partial_a^\omega=\partial_a-i\omega_a$.
From now on, we will denote $g_N$ simply as $g$ while the metric $g_{M_z}$ on $M_z$ will be denoted by $g_z$. Consistently, the metric $g_M$ will be denoted with $g_0$. Using this convention and defining $\omega_{ab}=\partial_a\omega_b-\partial_b\omega_a$ the operator $d^*_Nd_N$ is given by 
\begin{equation*}
	d_N^*d_N\omega=-\frac{1}{\sqrt{\det g_z}}g_{ab}\partial_{c}\left(\sqrt{\det g_z}g^{dc}g^{eb}\omega_{de}\right)dx^{a}.
\end{equation*}
From \eqref{gN} follows that 
\begin{equation*}
	g^{ij}=g_z^{ij},\quad g^{j\alpha}=0,\quad g^{\alpha\beta}=\delta_{\alpha\beta}.
\end{equation*}
We define
\begin{align*}
	a^{ij}&=g^{ij},\\ b_{s}^{ik}&=\frac{1}{\sqrt{\det g}}g_{st}\partial_{j}\left(g^{ij}g^{kt}\sqrt{\det g}\right),\\
	c^{k}&=\frac{1}{\sqrt{\det g}}\partial_{j}\left(\sqrt{\det g}g^{jk}\right),\\ d_{j}^{\beta k}&=\frac{1}{\sqrt{\det g}}g_{ij}\partial_{\beta}\left(g^{ik}\sqrt{\det g}\right),\\
	H_z^\beta& =-\frac{1}{\sqrt{\det g}}\partial_{\beta}\left(\sqrt{\det g}\right)
\end{align*}
and observe that $\sqrt{\det g}=\sqrt{\det g_z}$. Now we consider local coordinates on $M$ describing a neighbourhood of a generic $p\in M$, given by 
\begin{equation*}
	Y_p: B(0,\varrho)\subset\R^n\to M,\quad Y_p:\xi\to Y_p(\xi)
\end{equation*}
for some $\varrho>0$.
Denoting $\partial_{ab}^\omega=\partial_a^\omega\partial_b^\omega$, we obtain the following local expressions
\begin{align}
	-\Delta_N^\omega\phi&=-a^{ij}\partial^\omega_{ij}\phi-c^j\partial^\omega_j\phi-\partial^\omega_{\alpha\alpha}\phi+H_z^\alpha\partial^\omega_\alpha\phi\label{expansion delta}\\
	\begin{split}
		d_N^{*}d_N\omega&=-a^{ij}\partial_{j}\omega_{ik}d\xi^{k}-b_{k}^{ij}\omega_{ij}d\xi^{k}-\partial_{\beta}\omega_{\beta\gamma}dz^{\gamma}+H_{{z}}^{\beta}\omega_{\beta\gamma}dz^{\gamma}\label{expansion d*d}\\
	&\ -a^{ij}\partial_{j}\omega_{i\gamma}dz^{\gamma}-c^{i}\omega_{i\gamma}dz^{\gamma}-\partial_{\beta}\omega_{\beta k}d\xi^{k}-d_{k}^{\beta j}\omega_{\beta j}d\xi^{k}
	\end{split}
\end{align}
where all coefficients are evaluated at $(\xi,z)$, leaving implicit the composition with $Y_p$. We remark that $H^\beta_z$ is the mean curvature of $M_z$ in the direction of $\nu_\beta$ which can be expanded, using the minimality of $M$, as 
\begin{align}\label{mean curvature expansion}
\begin{split}
	H_z^\beta(\xi)&=z^{\gamma}\left(k_{\ell\beta}(\xi)k_{\ell\gamma}(\xi)+\Ric(\nu_{\beta},\nu_{\gamma})(\xi)\right)\\&\ +z^{\delta}z^{\gamma}\left(k_{\ell\beta}(\xi)k_{\ell\gamma}(\xi)k_{\ell\delta}(\xi)+\nabla\text{Ric}(\nu_{\delta},\nu_{\gamma};\nu_{\beta})(\xi)\right)+O(|z|^{3})
\end{split}
\end{align}
where $\Ric$ is the Ricci tensor of $N$ and $\nabla\Ric$ is its covatiant derivative. Here, $k_{\ell\gamma}$ is the $\ell$-th principal curvature in the direction of $\nu_\gamma$, precisely the $\ell$-th eigenvalue of the tensor $A_\gamma$ defined by 
\begin{equation*}
	A_\gamma(v,w)=-g(\nabla_v\nu_\gamma,w),\quad v,w\in TM,
\end{equation*}
see for instance \cite[Appendix A.2]{Badran-delPino2022}.
We write the expansions
\begin{align*}
	a^{ij}(\xi,z)&=a^{ij}_{0}(\xi)+z^\beta a^{ij}_{1,\beta}(\xi,z)\\
	c^j(\xi,z)&=c^j_{0}(\xi)+z^\beta c^j_{1,\beta}(\xi,z)
\end{align*}
and obtain
\begin{align*}
	a^{ij}\partial_{ij}+c^j\partial_j&=\Delta_M+\eps(t^\beta+h^\beta)\round{a^{ij}_{1,\beta}\partial_{ij}+c^j_{1,\beta}\partial_j}\nonumber\\
	&=:\Delta_M+\eps(t^\beta+h^\beta)D_{1,\beta}.\label{expansion laplacian}
\end{align*}
At this point we change coordinates to $X_h$. In what follows we denote $h_j^\beta=\partial_jh^\beta$, $h_{ij}^\beta=\partial_{ij}h^\beta$ and so on. Changing coordinates yields, for a pair $(\phi,\omega)$ locally represented as $\phi=\phi(\xi,t)$ and $\omega=\omega_{k}(\xi,t)d\xi^k+\omega_\alpha(\xi,t)dt^\alpha$, to 
\begin{equation}\label{expansion laplacian}
	\begin{split}
	-\eps^2\Delta^{\omega}_{N}\phi&=-\eps^2\Delta_{M_{z}}^{\omega}\phi-\partial_{\alpha\alpha}^{\omega}\phi+\eps H_{{z}}^{\beta}\partial_{\beta}^{\omega}\phi+\eps^{2}\left(\Delta_{M_{z}}h^{\beta}\right)\partial_{\beta}^{\omega}\phi\\&\ +\eps a^{jk}\left[h_{k}^{\beta}\partial_{j\beta}^{\omega}\phi+h_{j}^{\gamma}\partial_{\gamma k}^{\omega}\phi\right]-\eps^{2}a^{jk}h_{j}^{\gamma}h_{k}^{\beta}\partial_{\gamma\beta}^{\omega}\phi
\end{split}
\end{equation}
where we stress that $\Delta_{M_z}^\omega=a^{ij}(\xi,z)\partial^\omega_{ij}+c^j(\xi,z)\partial^\omega_j$ is a differential operator only in the variables $\xi$ and we recall that $\partial_j^\omega=\partial_j-i\omega_j$. It holds
\begin{equation}\label{expansion d*d}
	\begin{split}
	\eps^2d^{*}_Nd_{N}\omega&=-a^{ij}\Big(\eps^2\partial_{j}\omega_{ik}-\varepsilon h_{j}^{\beta}\partial_{\beta}\omega_{ik}-\varepsilon h_{i}^{\beta}\partial_{j}\omega_{\beta k}-\varepsilon^{2}h_{ji}^{\beta}\omega_{\beta k}\\&\ -\varepsilon^{2}h_{jk}^{\beta}\omega_{i\beta}+\varepsilon^{3}h_{i}^{\gamma}h_{jk}^{\beta}\omega_{\gamma\beta}+\varepsilon^{3}h_{i}^{\beta}h_{j}^{\gamma}\partial_{\gamma}\omega_{\beta k}\Big)d\xi^{k}\\&\ -\varepsilon b_{k}^{ij}\left(\eps^2\omega_{ij}-\varepsilon h_{i}^{\beta}\omega_{\beta j}-\varepsilon h_{j}^{\beta}\omega_{i\beta}+\varepsilon^{2}h_{i}^{\gamma}h_{j}^{\beta}\omega_{\gamma\beta}\right)d\xi^{k}-\partial_{\beta}\omega_{\beta k}d\xi^{k}\\&\ -\varepsilon d_{k}^{\beta j}\left(\eps\omega_{\beta j}-\varepsilon h_{j}^{\gamma}\omega_{\beta\gamma}\right)d\xi^{k}+\varepsilon^{2}h_{k}^{\gamma}H_z^{\beta}\omega_{\beta\gamma}d\xi^{k}-\varepsilon^{2}h_{k}^{\gamma}c^{i}\left(\omega_{i\gamma}-\varepsilon h_{i}^{\beta}\omega_{\beta\gamma}\right)d\xi^{k}\\&\ -\partial_{\beta}\omega_{\beta\gamma}dt^{\gamma}+\varepsilon H_{z}^{\beta}\omega_{\beta\gamma}dt^{\gamma}-\varepsilon c^{i}\left(\omega_{i\gamma}-\varepsilon h_{i}^{\beta}\omega_{\beta\gamma}\right)dt^{\gamma}\\&\ -a^{ij}\left(\eps\partial_{j}\omega_{i\gamma}-\varepsilon h_{j}^{\beta}\partial_{\beta}\omega_{i\gamma}-\varepsilon^{2}h_{ij}^{\beta}\omega_{\beta\gamma}-\varepsilon h_{i}^{\beta}\partial_{j}\omega_{\beta\gamma}+\varepsilon^{2}h_{i}^{\beta}h_{j}^{\delta}\partial_{\delta}\omega_{\beta\gamma}\right)dt^{\gamma}
\end{split}
\end{equation}
where all coefficients above have to be evaluated either at $(\xi)$ or $(\xi,t+h)$. A useful remark is that if the 1-form $\omega$ is of the form $\omega=\omega_\alpha(t)dt^\alpha$, then only the terms of the form $\omega_{\beta\gamma}$ don't vanish and hence
\begin{align*}
	\eps^2d^{*}_Nd_{N}\omega&=\left(-\partial_{\beta}\omega_{\beta\gamma}+\varepsilon H_{{z}}^{\beta}\omega_{\beta\gamma}+\varepsilon^{2}\left(\Delta_{M_{z}}h^{\beta}\right)\omega_{\beta\gamma}-\varepsilon^{2}a^{ij}h_{i}^{\beta}h_{j}^{\delta}\partial_{\delta}\omega_{\beta\gamma}\right)dt^{\gamma}\\&\ +\Big(-\varepsilon^{3}a^{ij}h_{i}^{\gamma}h_{jk}^{\beta}\omega_{\gamma\beta}-\varepsilon^{3}b_{k}^{ij}h_{i}^{\gamma}h_{j}^{\beta}\omega_{\gamma\beta}\\&\ +\varepsilon^{3}d_{k}^{\beta j}h_{j}^{\gamma}\omega_{\beta\gamma}+\varepsilon^{2}h_{k}^{\gamma}H_z^{\beta}\omega_{\beta\gamma}+\varepsilon^{3}h_{k}^{\gamma}c^{i}h_{i}^{\beta}\omega_{\beta\gamma}\Big)d\xi^{k}.
\end{align*}
With these formulas we can evaluate $S_\eps(W_0)$. Using the fact that the building block \eqref{building-block} solves the equations in $\R^2$, we obtain 
\begin{equation}\label{unexpanded first error}
	\begin{split}
	S_\eps(W_0)&=\eps^2(\Delta_{M_{z}}h^{\beta})(\xi)\sdf{V}_\beta(t)+\eps H^\beta_z(\xi) \sdf{V}_\beta(t)\\
	&\ -\eps^2a^{ij}(\xi,\eps(t+h))h_i^\beta(\xi) h_j^\gamma(\xi)\nabla_{\beta\gamma,U_0}U_0(t)+\eps^3R_\xi(h)
\end{split}
\end{equation}
where $\sdf{V}_\beta$ is as in \eqref{expression sdfVj} and $\nabla_{U_0}$ is given by \eqref{gradient-like}. Here, $R_\xi(h)$ is a remainder given by
\begin{equation*}
	R_\xi(h)=-\round{0,\round{a^{ij}h_{i}^{\gamma}h_{jk}^{\beta}+b_{k}^{ij}h_{i}^{\gamma}h_{j}^{\beta}+d_{k}^{\beta j}h_{j}^{\gamma}+\varepsilon^{-1}h_{k}^{\gamma}H_z^{\beta}+h_{k}^{\gamma}c^{i}h_{i}^{\beta}}a'd\xi^{k}}^T.
\end{equation*}
Next, we expand \eqref{unexpanded first error} around $M$. Using \eqref{mean curvature expansion} and \eqref{expansion laplacian}, we find 
\begin{equation}\label{SW0}
	\begin{split}
		S_\eps(W_0)&=\eps^2(\Delta_{M}h^{\beta})(\xi)\sdf{V}_\beta(t)+\eps^3(t^\gamma+h^\gamma)(D_{1,\gamma}h^{\beta})(\xi,\eps(t+h))\sdf{V}_\beta(t)\\
	&\ +\eps^2 (t^{\gamma}+h^\gamma)\left(k_{\ell\beta}(\xi)k_{\ell\gamma}(\xi)+\Ric(\nu_{\beta},\nu_{\gamma})(\xi)\right)\sdf{V}_\beta(t)\\
	&\ +\eps^3t^\gamma t^\delta k_{\ell\beta}(\xi)k_{\ell\gamma}(\xi)k_{\ell\delta}(\xi)\sdf{V}_\beta(t) + \eps^3t^\gamma t^\delta\nabla\Ric(\nu_\gamma,\nu_\beta;\nu_\delta)(\xi)\sdf{V}_\beta(t)\\
	&-\eps^2a_0^{ij}(\xi)h_i^\beta (\xi)h_j^\gamma(\xi)\nabla_{\beta\gamma,{U_0}}{U_0}(t)\\
	&\ -\eps^3(t^\delta+h^\delta)a_{1,\delta}^{ij}(\xi,\eps(t+h))h_i^\beta (\xi)h_j^\gamma(\xi)\nabla_{\beta\gamma,{U_0}}{U_0}(t)+\eps^3R_\xi(h)+O(\eps^{4}).
	\end{split}
\end{equation}
At this point we can try to improve the approximation by setting $W_1=W_0+\Lambda_1$ and observing that 
\begin{align*}
	S_\eps(W_1)=S_\eps(W_0)+\sdf{L}[\Lambda_1]+(S_\eps'(W_0)-\sdf{L})[\Lambda_1]+N_0(\Lambda_1)
\end{align*}
where 
\begin{equation*}
	N_0(\Lambda_1)=S_\eps(W_0+\Lambda_1)-S_\eps(W_0)-S_\eps'(W_0)[\Lambda_1]
\end{equation*}
where we recall that the 2-dimensional linearized operator $\sdf{L}$, given by \eqref{lstr}, is an operator in the $t$-variable only. The largest term of $S_\eps(W_0)$, namely that of order $\eps^2$, is locally given by
\begin{equation}\label{Q2}
	\sdf{Q}_2(\xi,t)=\eps^2\left(k_{\ell\beta}(\xi)k_{\ell\gamma}(\xi)+\Ric(\nu_{\beta},\nu_{\gamma})(\xi)\right)t^\gamma\sdf{V}_\beta(t).
\end{equation}
Therefore, if we solve 
\begin{equation*}
	\sdf{L}[\Lambda_1]=-\sdf{Q}_2(y,t)
\end{equation*}
the biggest part of the error in terms of $\eps$ is cancelled. Such $\Lambda_1$ exists thanks to Lemma \ref{invertibility lstr}, using that 
\begin{equation}\label{ort cond Q2}
	\int_{\R^2}\sdf{Q_2}(y,t)\cdot\sdf{V}_\alpha(t)dt=0, \quad \forall y\in M
\end{equation}
which follows directly from the fact that $\int_{\R^2}t^\gamma\sdf{V}_\gamma(t)\cdot\sdf{V}_\alpha(t)dt=0$. Also, since the right hands side \eqref{Q2} is $O(e^{-|t|})$ for $|t|$ large a standard barrier argument
along with the fact that $\sdf{L}\sim -\Delta + \mathrm{Id}$ at infinity ensures that
\begin{equation*}
	\sup_{t\in\R^2}e^{\sigma|t|}|\sdf{Q}_2(y,t)|<\infty\quad \forall y\in M
\end{equation*}
for any $\sigma<1$. Lastly, we observe that the error created 
\begin{equation*}
	\mathcal{E}(y,t)=(S_\eps'(W_0)-\sdf{L})[\Lambda_1]+N_0(\Lambda_1)
\end{equation*}
satisfies
\begin{equation}\label{error estimate}
	|\mathcal{E}(y,t)|\leq C\eps^4e^{-\sigma|t|}.
\end{equation}
Indeed, it holds an equation of the form 
\begin{equation*}
	(S_\eps'(W_0)-\sdf{L})[\Lambda_1]+N_0(\Lambda_1)=\alpha^{ij}D_{ij}\Lambda_1+\beta^jD_j\Lambda_1+O(\eps^4)
\end{equation*}
with $\alpha^{ij},\beta^j=O(\eps^2)$ and hence, using the exponential decay in $|t|$ and the fact that $\Lambda_1=O(\eps^2)$, we check the validity of \eqref{error estimate}.
Thus, with the above improvement for the first approximation the following estimate holds
\begin{equation*}
	\|S_\eps(W_1)\|_{C^{0,\gamma}}\leq C \eps^3.
\end{equation*}

We can obtain a further improvement of the approximation by eliminating the $\eps^3$-terms in $S_\eps(W_1)$.
To do this, the non-degeneracy assumption \eqref{nondeg} becomes crucial. 
Proceeding as before, we look for a bounded function $\Lambda_2$ such that
\begin{equation}\label{L Lambda1}
	\sdf{L}[\Lambda_2] =  -\sdf{Q}_3(y,t)
\end{equation}
where $\sdf{Q}_3$ are the $\eps^3$-terms in $S(W_1)$, which from \eqref{SW0} and \eqref{assumption h} are
\begin{align*}
	\sdf{Q}_3(y,t)&=\eps^2\round{\Delta_Mh^\beta+k_{\ell\gamma}k_{\ell\beta}h^\gamma+\Ric(\nu_\beta,\nu_\gamma)h^\gamma}\sdf{V}_\beta(t)\\
	&\ +\eps^3\round{k_{\ell\beta}(y)k_{\ell\gamma}(y)k_{\ell\delta}(y) + \nabla\Ric(\nu_\gamma,\nu_\beta;\nu_\delta)(y)}t^\gamma t^\delta\sdf{V}_\beta(t)
\end{align*}
and then set $W_2=W_1+\Lambda_2$. This can be done by virtue of Lemma \ref{invertibility lstr} if the right-hand side satisfies the orthogonality condition \eqref{lemma1 ort cond}, namely if the quatities
\begin{align}\label{inner product rhs}
	\int_{\R^2}\sdf{Q}_3(y,t)\cdot \sdf{V}_\gamma(t)\,  dt = c\eps^2\J^\gamma(h)(y)+\eps^3\sdf{q}^\gamma(y),
\end{align}
vanish for $\gamma=1,2$, where
\begin{equation*}
	\J^\gamma[h]=\Delta_Mh^\gamma+k_{\ell\beta}k_{\ell\gamma}h^\beta+\Ric(\nu_\gamma,\nu_\beta)h^\beta,\quad \gamma=1,2
\end{equation*}
are the components of the Jacobi operator, $c=\int_{\R^2}|\sdf{V}_\gamma(t)|^2dt$ (independent of $\gamma=1,2$) and 
\begin{equation}\label{qgamma}
	\sdf{q}^\gamma(y) = \round{k_{\ell\gamma}(y)k_{\ell\beta}(y)k_{\ell\delta}(y) + \nabla\Ric(\nu_\beta,\nu_\gamma;\nu_\delta)(y)}\int_{\R^2}t^\delta t^\beta|\sdf{V}_\gamma(t)|^2dt
\end{equation}
Assumption \eqref{nondeg} of non-degeneracy and Fredholm alternative for elliptic operators imply the following result.
\begin{lemma}\label{invert Jacobi}
	Let $f\in C^{0,\gamma}(M,\R^2)$. Then the system
	\begin{equation*}
		\J(h)=f\quad\text{on }M
	\end{equation*}
	admits a solution $h=\mathcal{H}(f)$ satisfying 
	\begin{equation*}
		\|h\|_{C^{2,\gamma}(M)}\leq C\|f\|_{C^{0,\gamma}(M)}.
	\end{equation*}
\end{lemma}
Lemma \ref{invert Jacobi} guarantees the existence of a bounded  $h_0(y)$ with 
\begin{equation}\label{eq h0}
    \J(h_0) = -(\sdf{q}^1,\sdf{q}^2)^T, \quad \text{on }M.
\end{equation}
Choosing  $h= \eps h_0$, the right hand side of 
\eqref{inner product rhs} vanishes for $\gamma=1,2$. Thus, 
\eqref{L Lambda1} has a unique solution $\Lambda_2 (y,t)$ with 
$$
\int_{\R^2} \Lambda_2 (y,t)\cdot \sdf{V}_\gamma (t)\, dt =0.
$$
The role of the invertibility on the Jacobi operator $\J$ is fairly transparent in the above argument, but since the operator $\sdf{L}$
does not include the differential operator in the $y$-variable, a linear theory for the full linearized  problem needs to be made. This is the content of Proposition 
\ref{inver L_W} below. The solvability conditions needed on the right hand sides only involve the $t$-variable and read very much as \eqref{lemma1 ort cond} or \eqref{ort cond Q2}.    
The approximation $W_1(y,t)$ will actually be sufficient for our purposes. The function, $\Lambda_2(y,t)$ will be part of the expansion of the full solution. 

The approximation obtained is only local and we need to extend it to the entire ambient manifold $N$. This extension creates very small errors because of the exponential decays involved. Its existence uses in essential way the orientability of $M$. 
Let $\delta>0$ and consider the cut-off functions defined by 
\begin{equation}\label{cutoffs}
	\zeta_m(x)=\begin{cases}
		\zeta(\frac{\eps}{\delta}|t+h(y)|-m) & \text{if }x=X_h(y,t)\in \N_h\\
		0 & \text{otherwise.}
	\end{cases}
\end{equation}
where $\zeta$ is a smooth cut-off function such that $\zeta(s)=1$ if $s<1$ and $\zeta(s)=0$ if $s>2$. 
To obtain a global approximation we first use the following Lemma, whose proof is the content of Appendix \ref{A1}.
\begin{lemma}\label{existence psi}
The admissibility hypothesis (H) guarantees the existence of a smooth function 
\begin{equation*}
	\psi : N\setminus M_h\to S^1,
\end{equation*}
where $M_h=\{\exp_y(h^\beta(y)\nu_\beta(y))\,:\,y\in M\}$,
 such that for every $x=X_h(y,t)\in \mathrm{supp}\,\zeta_3$
 \begin{equation*}
 	\psi(x) =\frac{t}{|t|},\quad \ |t+h(y)|< 5\delta/\eps.
 \end{equation*}
\end{lemma}
We consider the pure gauge
\begin{equation*}
	\boldsymbol{\Psi}=\twovec{\psi}{\frac{d\psi}{i\psi}}
\end{equation*}
where $\psi$ is the smooth function in Lemma \ref{existence psi} and define the global approximation $W$ to a solution of \eqref{equations} as
\begin{equation*}
	W=\zeta_3W_1+(1-\zeta_3)\boldsymbol{\Psi}.
\end{equation*} 
A lengthy but straightforward calculation similar to that of \cite[\S3.7]{Badran-delPino2022} shows that, for $\sigma\in(0,1)$
\begin{equation*}
	S_\eps(W)=\zeta_3S_\eps(W_1)+(1-\zeta_3)S_\eps(\boldsymbol{\Psi})+\mathrm{E}
\end{equation*}
where $S_\eps(\boldsymbol{\Psi})=0$ since $\boldsymbol{\Psi}$ is a pure gauge and 
\begin{equation*}
	|\mathrm{E}|\leq C\chi_{\{0<\zeta_3<1\}}(1+|W_1|^2)(|W_1-\boldsymbol{\Psi}|+|D(W_1-\boldsymbol{\Psi})|+|\nabla^{A_1}u_1|).
\end{equation*}
Now, by construction of $\Psi$ and the properties of the function $f$ and $a$ introduced in \eqref{building-block}, we have that
\begin{equation*}
	|W_1-\boldsymbol{\Psi}|+|D(W_1-\boldsymbol{\Psi})|+|\nabla^{A_1}u_1|\leq Ce^{-\sigma|t|}\quad \text{on }\{0<\zeta_3<1\}.
\end{equation*}
Recalling that on such set we have 
\begin{equation*}
	e^{-\sigma|t|}\leq e^{-\frac{4\sigma\delta}{\eps}}
\end{equation*}
we obtain that $|\mathrm{E}|$ is exponentially small in $\eps$.  
\section{Proof of main result}
Starting from the global approximation just constructed we look for an actual solution as a small perturbation of $W$, namely we want to solve
\begin{equation}\label{perturbed solution}
	S_\eps(W+\vphi)=0
\end{equation}
for some $\vphi$ which is small in a suitable sense. In \S\ref{role gauge} we proved that to solve \eqref{perturbed solution} it suffices to solve 
\begin{equation}\label{proj perturbed solution}
	\pi_{Z_{W,g}^\perp}S_\eps(W+\vphi)=0
\end{equation}
for a $\vphi\in Z_{W,g}^\perp$ which is sufficiently small with respect to $W$. Using the decomposition \eqref{decomposition} of the linearized and the fact that $S_\eps(W)\in Z_{W,g}^\perp$ we find that \eqref{proj perturbed solution} is equivalent to 
\begin{equation}\label{uncorrected}
	L^\eps_W[\vphi]=-S_\eps(W)-\pi_{Z_{W,g}^\perp}N_\eps(\vphi)\quad\text{if }\vphi\in Z_{W,g}^\perp
\end{equation}
where 
\begin{equation*}
	N_\eps(\vphi)\coloneqq S_\eps(W+\vphi)-S_\eps(W)-S_\eps'(W)[\vphi]
\end{equation*}
and $L^\eps_W$ is given by \eqref{LW}. More precisely,
\begin{equation*}
	L^\eps_W[\vphi]=-\eps^2\Delta_W\vphi+\vphi+T_W^\eps\vphi
\end{equation*}
Rather than \eqref{uncorrected}, we solve a corrected version, given by
\begin{equation}\label{corrected}
	L^\eps_W[\vphi]=-S_\eps(W)-\pi_{Z_{W,g}^\perp}N_\eps(\vphi)+\zeta_2b^\alpha(y)\sdf{V}_\alpha(t),\qquad\vphi\in Z_{W,g}^\perp.
\end{equation}
The adjustment on the right-hand side provides unique solvability in terms of $\vphi$ for a precise choice of $b=(b^1,b^2)$, in the sense of the following result.
\begin{proposition}\label{inver L_W}
	Let $0<\gamma<1$ and let $\Lambda\in C^{0,\gamma}(N)$. Then there exists $b\in C^{0,\gamma}(M)$ and a unique solution $\vphi=\mathcal{G}(\Lambda)$ to
	\begin{equation*}
		L^\eps_W[\vphi]=\Lambda+\zeta_2b^\alpha(y)\sdf{V}_\alpha(t)
	\end{equation*}
	satisfying
	\begin{equation*}
		\|\vphi\|_{C^{2,\gamma}(N)}+\|b\|_{C^{0,\gamma}(M)}\leq C \|\Lambda\|_{C^{0,\gamma}(N)}
	\end{equation*}
\end{proposition}

Proposition \ref{inver L_W} allows us to write \eqref{corrected} as a fixed point problem 
\begin{equation}\label{fixed point formulation}
	\vphi=-\mathcal{G}\round{S_\eps(W)+\pi_{Z_{W,g}^\perp}N_\eps(\vphi)}
\end{equation}
on the space 
\begin{equation*}
	X_A=\braces{\vphi\in C^{2,\gamma}(N)\ :\ \|\vphi\|_{C^{2,\gamma}(N)}\leq A\eps^3}\cap Z_{W,g}^\perp.
\end{equation*}
Using Remark \eqref{self-map} and the Lipschitz character of $N_\eps$ on $X_A$, namely
\begin{equation*}
	\|N_\eps(\vphi_1)-N_\eps(\vphi_2)\|_{C^{2,\gamma}(N)}\leq C\eps^3\|\vphi_1-\vphi_2\|_{C^{2,\gamma}(N)},\quad \vphi_1,\vphi_2\in X_A,
\end{equation*}
we find a solution $\vphi\in X_A$ of \eqref{fixed point formulation}. We also highlight the Lipschitz dependence of $\vphi$ on $h$, namely 
\begin{equation*}
	\|\vphi(h_1)-\vphi(h_2)\|_{C^{2,\gamma}(N)}\leq C\eps^2\|h_1-h_2\|_{C^{2,\gamma}(M)}
\end{equation*}
which will be proved in \S\ref{G contract section}. The only step left is to adjust $h$ in a way that the correction $\zeta_2b^\alpha\sdf{V}_\alpha$ vanishes. 
To do so, we multiply \eqref{corrected} by $\zeta_4\sdf{V}_\gamma(t)$, $\gamma=1,2$, and integrate on $\R^2$ to find an expression for $b^\gamma$,
\begin{equation*}
	b^\gamma(y)=\frac{1}{\int_{\R^2}\zeta_2|\sdf{V}_\gamma|^2}\int_{\R^2}\squared{\zeta_4S_\eps(W)+\zeta_4\pi_{Z_{W,g}^\perp}N(\vphi)+\zeta_4L_W[\vphi]}\cdot\sdf{V}_\gamma.
\end{equation*}
Working out expression \eqref{SW0} we find that, also thanks to assumption \eqref{assumption h}, the biggest term of the error $S_\eps(W)$ within the support of $\zeta_4$ is $\eps^2 \J^\gamma[h]\sdf{V}_\gamma$, where
\begin{equation}
	\J^\gamma[h]=\Delta_Mh^\gamma+k_{\ell\beta}k_{\ell\gamma}h^\beta+\Ric(\nu_\gamma,\nu_\beta)h^\beta,\quad \gamma=1,2
\end{equation}
are the components of the Jacobi operator. Let us set 
\begin{equation*}
	q_m(y)=\int_{\R^2}\zeta_m(y,t)|\sdf{V}_\alpha(t)|^2,\quad m=1,2,\dots
\end{equation*}
and observe that the above quantity does not depend on $\alpha=1,2$. We define the nonlinear operator $G$ in the following way: for $\alpha=1,2$, let
\begin{align*}
	G^1_\alpha(h)&=q_4\J^\alpha[h]-\eps^{-2}\int_{\R^2}\squared{\zeta_4S(W_1)}\cdot\sdf{V}_\alpha,\\
	G^2_\alpha(h)&=-\eps^{-2}\int_{\R^2}\zeta_4\squared{\pi_{Z_{W,g}^\perp}N(\vphi)+\sdf{B}[\vphi]}\cdot\sdf{V}_\alpha,\\
	G^3_\alpha(h)&=-\eps^{-2}\int_{\R^2}\zeta_4L_{U_0}[\vphi]\cdot\sdf{V}_\alpha
\end{align*}
and set
\begin{equation*}
	G(h)=\twovec{G_1(h)}{G_2(h)}=q_4^{-1}\sum_{k=1}^3\twovec{G_1^k(h)}{G_2^k(h)}.
\end{equation*}
With this notation, the system $b^\gamma=0$, $\gamma=1,2$, reads
\begin{equation}\label{Jacobi to solve}
	\J(h)=G(h).
\end{equation}
Observe that $G^1$ satisfies
\begin{equation*}
	G^1_\alpha(y)=- \round{k_{\ell\alpha}(y)k_{\ell\beta}(y)k_{\ell\delta}(y) + \nabla\Ric(\nu_\beta,\nu_\alpha;\nu_\delta)(y)}\int_{\R^2}\zeta_4(y,t)t^\delta t^\beta|\sdf{V}_\alpha(t)|^2dt+O(\eps)
\end{equation*} 
and its main order corresponds to the negative of the term \eqref{qgamma} previously described (with a slight variation due to the introduction of the cut-off $\zeta_4$).
To solve \eqref{Jacobi to solve} we use Lemma \ref{invert Jacobi}  and restate it as a fixed point problem 
\begin{equation}\label{fixed point h}
	h=\mathcal{H}(G(h)).
\end{equation}
To conclude the proof we use the following Lemma, which will be proved in \S\ref{G contract section}.
\begin{lemma}\label{G contracts}
	The map $G$ satisfies 
	\begin{equation*}
		\|G(0)\|_{C^{0,\gamma}(M)}\leq C\eps
	\end{equation*}
	and 
	\begin{equation*}
		\|G(h_1)-G(h_2)\|_{C^{0,\gamma}(M)}\leq C\eps\|h_1-h_2\|_{C^{2,\gamma}(M)},
	\end{equation*}
	for some $C>0$.
\end{lemma}
Thanks to Lemma \ref{G contracts} we see that by contraction mapping principle equation \eqref{fixed point h} admits a solution in the space 
\begin{equation*}
	X_A=\braces{h\ :\ \|h\|_{C^{2,\gamma}(M)}\leq A\eps},
\end{equation*}
for any $A$ sufficiently big, which concludes the proof. Let us observe that, looking at the expansion of the terms of $G(h)$ in \eqref{Jacobi to solve}, we have $h(y)=\eps h_0(y)+O(\eps^2)$ where $h_0(y)$ is the unique solution of equation \eqref{qgamma}-\eqref{eq h0}.

\medskip
Finally, we consider the zero-set of $u_\eps(x)$. The function $u_\eps$ does not vanish away from $M$. Using \eqref{local info}, the equation $u_\eps (x)=0$ becomes near $M$ 
$$
u_0( \eps^{-1}(z- \eps^2h_0(y)))  + \eps^2\theta(\eps^{-1}z,y)  =0   
$$
where $h_0(y)$ is the unique solution of 
and $\theta$ is smooth, with uniformly bounded derivatives. The fact that $\partial_zu_0(0)\ne0$  and the implicit function theorem
yield that the zero set is a smooth manifold parametrized as
$$
z= \eps^2h_0(y)+  O(\eps^3)  , \quad y\in M.
$$
The proof is complete. 
\qed
\section{Proof of Proposition \ref{inver L_W}}
To prove Proposition \ref{inver L_W} we will use the fact that on a region close to $M$ the linearized $L_W$ can be approximated by $L_{U_0}$, namely the scaled linearized operator on $M\times\R^2$ around the building block $U_0(y,t) \coloneqq U_0(t)$, which in the scaled coordinates $(y,t)=(y,z/\eps)$ is given by
\begin{equation*}
	L^\eps_{U_0}[\Phi]=-\Delta_{t,U_0}\Phi-\eps^2\Delta_{M}\Phi+\Phi+T_{U_0}(t)\Phi,
\end{equation*}
being
\begin{equation*}
	T_{U_0}=\twovec{-\frac12(1-|u_0|^2)\phi+2i\nabla^{A_0}u_0\cdot\omega}{-(1-|u_0|^2)\omega+2\langle\nabla^{A_0}u_0,i\phi\rangle}.
\end{equation*}
First, we consider the cut-off functions introduced in \eqref{cutoffs} and let us look for a solution to 
\begin{equation}\label{CorLinePro1}
	L^\eps_W[\vphi]=-\eps^2\Delta_W\vphi+\vphi+T^\eps_W\vphi=\Lambda+\zeta_2b^\alpha\sdf{V}_\alpha\quad \text{on }N
\end{equation}
of the form
\begin{equation*}
	\vphi(x)=\zeta_2(x)\Phi(y,t)+\Psi(x)
\end{equation*}
where $\Phi$ is defined on $M\times \R^2$ and $\Psi$ is defined on $N$. In what follows it will be useful the following definition: for a real valued function $f$ we set
\begin{align*}
	\mathcal{R}_W[f,\Phi]=-\Delta_W(f\Phi)+f\Delta_W\Phi.
\end{align*}
More precisely, if $\Phi=(\phi,\omega)^T$ and $W=(u,A)^T$, we have
\begin{equation}\label{definition of R}
	\mathcal{R}_W[f,\Phi]=\begin{pmatrix}-2\nabla^{A}\phi\cdot d f-\phi\Delta f\\
d^{*}\left(df\wedge\omega\right)+d^{*}\omega df-d\left(df\cdot\omega\right)-\omega\Delta f
\end{pmatrix}
\end{equation}
and it holds
\begin{equation*}
	\abs{\mathcal{R}_W[f,\Phi]}\leq C(|df|+|D^2 f|)(|\Phi|+|\nabla_W\Phi|).
\end{equation*} Equation \eqref{CorLinePro1} becomes
\begin{align*}
	L^\eps_{W}[\vphi]&=\zeta_2L^\eps_{U_0}[\Phi]+L^\eps_W[\Psi]+\zeta_2(L^\eps_W-L^\eps_{U_0})[\Phi]+\eps^2\mathcal{R}_W[\zeta_2,\Phi]\\
	&=\Lambda+\zeta_2b^\alpha\sdf{V}_\alpha
\end{align*}
and such equation is solved if the pair $(\Phi,\Psi)$ solves the system
\begin{align}
	&L^\eps_{U_0}[\Phi]+(L^\eps_W-L^\eps_{U_0})[\Phi]+\zeta_1T^\eps_W\Psi=\Lambda+b^\alpha\sdf{V}_\alpha\quad\text{on supp}\,\zeta_2,\label{first_eq}\\
	&-\eps^2\Delta_W\Psi+\Psi+(1-\zeta_1)T^\eps_W\Psi+\eps^2\mathcal{R}_W[\zeta_2,\Phi]=(1-\zeta_2)\Lambda\quad\text{on }N\label{second_eq},
\end{align}
where we used that $\zeta_2\zeta_1=\zeta_1$. The following Lemma,  will allow us to reduce the above system to an equation depending only on $\Phi$.

\begin{lemma}\label{Outer invert}
	For $\gamma\in(0,1)$ and every $\eps>0$ and $\Gamma\in C^{0,\gamma}(N)$ there exists only one $\Psi$ satisfying 
	\begin{equation*}
		-\eps^2\Delta_W\Psi+\Psi+(1-\zeta_1)T^\eps_W\Psi=\Gamma\quad\text{on }N
	\end{equation*}
	and such that 
	\begin{equation*}
		\|\Psi\|_{C^{2,\gamma}(N)}\leq C\|\Gamma\|_{C^{0,\gamma}(N)}.
	\end{equation*}
\end{lemma}
Lemma \ref{Outer invert} follows from the positivity of the operator on the left hand side on the space $H^1_W(N)$ defined by the norm in \eqref{H1N}, the maximum principle and elliptic estimates. 
We define the weighted norm $C^{0,\gamma}_\sigma$ on functions $\psi(y,t)$ defined on $M\times\R^2$ as 
\begin{equation*}
	\|\psi\|_{C^{k,\gamma}_\sigma(M\times\R^2)}=\|e^{\sigma|t|}\psi\|_{C^{k,\gamma}(M\times\R^2)},
\end{equation*}
where $k\geq 0$ and $\gamma,\sigma\in(0,1)$.
Using the fact that on the support of $\zeta_2$, given also hypothesis \eqref{assumption h}, it holds $|t|>\delta/\eps$ we find that 
\begin{equation*}
	e^{-\sigma|t|}\leq e^{-\frac{\sigma\delta}{\eps}}\quad\text{on }\{0<\zeta_2<1\}
\end{equation*}
and hence
\begin{equation*}
	\|\mathcal{R}_W[\zeta_2,\Phi]\|_{C^{0,\gamma}(N)}\leq Ce^{-\frac{\sigma\delta}\eps}\|\Phi\|_{C_\sigma^{0,\gamma}(M\times\R^2)}.
\end{equation*}
Using Lemma \ref{Outer invert} we find a solution $\Psi=\Psi_1+\Psi_2$ to \eqref{second_eq}, where $\Psi_1,\Psi_2$ satisfy 
\begin{align*}
	-\eps^2\Delta_W\Psi_1+\Psi_1+(1-\zeta_1)T^\eps_W\Psi_1&=-\mathcal{R}_W[\zeta_2,\Phi],\\
	-\eps^2\Delta_W\Psi_2+\Psi_2+(1-\zeta_1)T^\eps_W\Psi_2&=(1-\zeta_2)\Lambda
\end{align*}
on $N$. Also, $\Psi$ satisfies 
\begin{equation*}
	\|\Psi\|_{C^{2,\gamma}(N)}\leq C\round{\|(1-\zeta_1)\Lambda\|_{C^{0,\gamma}(N)}+e^{-\frac{\sigma\delta}\eps}\|\Phi\|_{C_\sigma^{0,\gamma}(M\times\R^2)}}.
\end{equation*}
At this point we reduce the whole system \eqref{first_eq}-\eqref{second_eq} to a single equation on $M\times\R^2$. We define 
\begin{equation*}
	\tilde{\sdf{B}}[\Phi]=\zeta_4{\sdf{B}}[\Phi]=\zeta_4(L^\eps_W-L^\eps_{U_0})[\Phi],\quad\tilde\Lambda=\zeta_4\Lambda,\quad (y,t)\in M\times\R^2.
\end{equation*}
We claim that $\tilde{\sdf{B}}$ satisfies
\begin{equation}\label{estimate Btilde}
	\|\tilde{\sdf{B}}[\Phi]\|_{C^{0,\gamma}(M\times\R^2)}\leq C\delta\|\Phi\|_{C^{2,\gamma}(M\times\R^2)}
\end{equation}
where $\delta$ is the one from the definition of $\zeta_4$ in \eqref{cutoffs}. 
For instance, we readily find from \eqref{expansion laplacian} that 
\begin{equation*}
	|\eps^2\Delta^{\omega}_{N}\phi(y,t)-(\eps^2\Delta_{M}^{\omega}\phi(y,t)+\partial_{\alpha\alpha}^{\omega}\phi(y,t))|\leq C \delta\|\phi\|_{C^{2}(M\times\R^2)},
\end{equation*}
and similarly for the other terms. 
With this notation, the system is equivalent to 
\begin{equation}\label{mixed-eq}
	L_{U_0}[\Phi]+\tilde{\sdf{B}}[\Phi]+\zeta_1T_W\Psi=\tilde{\Lambda}+b^\alpha\sdf{V}_\alpha\quad \text{on } M\times\R^2.
\end{equation}
We will prove the following result.
\begin{proposition}\label{invert LU0}
	Let $\gamma\in (0,1)$ and $\sigma>0$ sufficiently small. Then for every $\tilde\Lambda\in C^{0,\gamma}_\sigma(M\times\R^2)$ there exists $b\in C^{0,\gamma}(M)$ such that the problem
\begin{equation*}
	L^\eps_{U_0}[\Phi]=\tilde\Lambda+b^\alpha\sdf{V}_\alpha\quad\text{on }M\times\R^2
\end{equation*}
admits a unique solution $\Phi=\mathcal{T}(\tilde\Lambda)$ satisfying 
\begin{equation*}
	\|\Phi\|_{C^{2,\gamma}_\sigma(M\times\R^2)}+\|b\|_{C^{0,\gamma}(M)}\leq C\|\tilde\Lambda\|_{C^{0,\gamma}_\sigma(M\times\R^2)}
\end{equation*}
for some $C>0$.
\end{proposition}
Using Proposition \ref{invert LU0} we can rephrase \eqref{mixed-eq} as 
\begin{equation}\label{phi+gphi}
	\Phi+\mathcal{G}[\Phi]=\mathcal{H}
\end{equation}
where
\begin{align*}
	\mathcal{G}[\Phi]&=\mathcal{T}\round{\tilde{\sdf{B}}[\Phi]+\zeta_1T^\eps_W\Psi_1[\Phi]}\\
	\mathcal{H}&=\mathcal{T}\round{\tilde\Lambda-\zeta_1T^\eps_W\Psi_2[\Lambda]}
\end{align*}
We observe that, thanks to \eqref{estimate Btilde}
\begin{equation*}
	\|\tilde{\sdf{B}}[\Phi]\|_{C^{0,\gamma}_{\sigma}(M\times\R^2)}\leq C\delta\|\Phi\|_{C^{2,\gamma}_{\sigma}(M\times\R^2)}
\end{equation*}	
where we obtain the control on exponential decay from the fact that $\tilde{\sdf{B}}[\Phi]$ is supported on $\mathrm{supp}\,\zeta_4$. Also 
\begin{equation*}
	\|\zeta_1T^\eps_W\Psi_1[\Phi]\|_{C^{0,\gamma}_{\sigma}(M\times\R^2)}\leq C \left\Vert \Psi_1[\Phi]\right\Vert _{C^{0,\gamma}(N)}\leq Ce^{-\frac{\delta'}{\eps}}\|\Phi\|_{C^{2,\gamma}_{\sigma}(M\times \R^2)}
\end{equation*}
where we used the fact that $\zeta_1T^\eps_W\sim e^{-|t|}\chi_{\{\zeta_1>0\}}(t)$ to control the exponential decay in the weighted norm. 
Thus, we have 
\begin{align*}
	\norm{\mathcal{G}[\Phi]}_{C^{2,\gamma}_{\sigma}(M\times\R^2)}&\leq C\round{\|\tilde{\sdf{B}}[\Phi]\|_{C^{0,\gamma}_{\sigma}(M\times\R^2)}+\|\zeta_1T_W\Psi_1[\Phi]\|_{C^{0,\gamma}_{\sigma}(M\times\R^2)}}\\
	&\leq C\round{\delta+e^{-\frac{\delta'}{\eps}}}\|\Phi\|_{C^{2,\gamma}_{\sigma}(M\times\R^2)}
\end{align*}
and hence, choosing $\eps$ and $\delta$ sufficiently small, we find a unique solution to \eqref{phi+gphi}, from which it follows the existence of a unique solution $(\Phi,\Psi)$ to system \eqref{first_eq}-\eqref{second_eq}. In conclusion, $\vphi=\zeta_2\Phi+\Psi$ solves \eqref{CorLinePro1} and  it follows directly that 
\begin{equation}
	\norm{\vphi}_{C^{2,\gamma}(N)}\leq C\norm{\Lambda}_{C^{0,\gamma}(N)}.
\end{equation}
which concludes the proof of Proposition \ref{inver L_W}.
\qed
\section{Invertibility of $L_{U_0}$}

In this section we prove a solvability theory for the equation 
\begin{equation}\label{uncorrected MR2}
	L^\eps_{U_0}[\Phi]=-\Delta_{t,U_0}\Phi-\eps^2\Delta_{M}\Phi+\Phi+T_{U_0}(t)\Phi=\Psi\quad\text{on }M\times\R^2.
\end{equation}
In general, it is not possible to solve \eqref{uncorrected MR2} for every choice of right-hand side $\Psi$. We solve instead the projected problem 
\begin{equation}\label{corrected MR2}
\begin{cases}
	L_{U_0}^\eps[\Phi]=\Psi+b^\alpha(y)\sdf{V}_\alpha(t)&\text{on }M\times\R^2\\
	\int_{\R^2}\Phi(y,t)\cdot\sdf{V}_\alpha(t)dt=0&\alpha=1,2,
\end{cases}	
\end{equation}
where 
\begin{equation}\label{balpha}
	b^\alpha(y)=-\frac{1}{\int_{\R^2}|\sdf{V}_\alpha(t)|^2dt}\int_{\R^2}\Psi(y,t)\cdot\sdf{V}_\alpha(t)dt,\quad y\in M,\ \alpha=1,2.
\end{equation}
This variation will provide unique solvability in the sense of Proposition \ref{invert LU0}. To prove this result we restrict to an open cover $\{\U_k\}$ of $M$ and solve the problem locally on $\U_k\times\R^2$ for every $k$, finding then a global solution by gluing of all the local solutions. First, we need an invertibility theory for the same operator on the flat space $\R^n$
\begin{equation}\label{corrected R4}
\begin{cases}
	-\Delta_{t,U_0}\Phi-\eps^2\Delta_{y}\Phi+\Phi+T_{U_0}\Phi=\Psi+b^\alpha(y)\sdf{V}_\alpha(t)&\text{on }\R^n\\
	\int_{\R^2}\Phi(y,t)\cdot\sdf{V}_\alpha(t)dt=0&\alpha=1,2,
\end{cases}	
\end{equation}
where $b^\alpha$ is given by \eqref{balpha}, $\alpha=1,2$ and $(y,t)\in\R^{n-2}\times\R^2$.
It holds the following Lemma, whose proof is postponed to Section \ref{proof of lemmas}.
\begin{lemma}\label{Solvability lin flat space}
	Let $\Psi\in L^2(\R^n)\cap C^{0,\gamma}(\R^n)$ and let $b^\alpha$ be given by \eqref{balpha} for $\alpha=1,2$. Then Problem \eqref{corrected R4}
admits a unique solution $\Phi=\mathcal{T}(\Psi)\in H^{1}_{U_0}(\R^n)$ satisfying 
\begin{equation}\label{holder estimates}
	\norm{\Phi}_{C^{2,\gamma}(\R^n)}\leq C\norm{\Psi}_{C^{0,\gamma}(\R^n)}
\end{equation}
for some $C>0$. 
\end{lemma}

\subsection{Proof of Proposition \ref{invert LU0}}. Given any point $p\in M$ we can find a local parametrization 
\begin{equation}\label{coordinates on M}
	Y_p:B(0,1)\subset\R^{n-2}\to M,\quad \xi \mapsto Y_p(\xi)
\end{equation}
such that 
\begin{equation*}
	g_{ij}(\xi)=\langle\partial_iY_p(\xi),\partial_jY_p(\xi)\rangle=\delta_{ij}+\theta_p(\xi),\quad \xi\in B(0,1)
\end{equation*}
where $\theta_p(0)=0$ and 
\begin{equation*}
	|D^2\theta_p|\leq C\quad \text{in }B(0,1). 
\end{equation*}
Locally, the Laplace-Beltrami operator reads 
\begin{equation}\label{expansion laplace-beltrami}
	\begin{split}
		\Delta_M&=\frac{1}{\sqrt{\det g}}\partial_i\round{\sqrt{\det g}g^ij\partial_j}\\
	&=\Delta_\xi+B_p
	\end{split}
\end{equation}
where $B_p$ has the form 
\begin{equation*}
	B_p=b^{ij}(\xi)\partial_{ij}+b^j(\xi)\partial_j,\quad |\xi|<1,
\end{equation*}
and $b^{ij},b^j$ and their derivatives are uniformly bounded. Let us fix $\delta>0$ sufficiently small and choose a sequence of points $(p_j)_{j\in\mathbb{N}}\subset M$ such that, if we define 
\begin{equation*}
	\V_j=Y_{p_j}(B(0,\delta/2)),\quad j=1,2,\dots
\end{equation*}
then $M$ is covered by the union of $\V_k$ and such that each $\V_k$ intersects at most a finite, uniform number of $\V_j$ with $j\ne k$. For a cut-off $\eta$ such that $\eta(s)=1$ if $s<1$ and $\eta(s)=0$ if $s>2$ define the following set of cut-offs on $M$
\begin{equation*}
	\eta_m^k(y)=\eta\round{\frac{|\xi|}{m\delta}},\quad y\in Y_{p_k}(\xi)
\end{equation*}
which are supported in $\U_{k,m}\coloneqq \braces{\xi\in\R^2\ :\ |\xi|\leq 2m\delta}$ and extended as $0$ outside of $\U_{k,m}$.
Observe that, given our choice of $\{\V_k\}$, there exists a $C>0$ for which
\begin{equation*}
	1\leq \eta_1\coloneqq\sum_{k=1}^\infty\eta_1^k\leq C.
\end{equation*}
We look for a solution to \eqref{corrected MR2} of the form 
\begin{equation*}
	\Phi=\Phi_0+\eta_1^k\Phi_k,\quad b^\alpha=\eta_1^kb_k^\alpha,\quad \alpha=1,2,
\end{equation*}
where we left implicit the sum over $k\in\mathbb{N}$ and the functions $\Phi_k,b^\alpha_k$ are defined in $\R^n$. With this ansatz, the equation can be written as 
\begin{align*}
	-\Delta_{t,U_0}\Phi-\eps^2\Delta_{M}\Phi+\Phi+T_{U_0}(t)\Phi&=-\Delta_{t,U_0}\Phi_0-\eps^2\Delta_{M}\Phi_0+\Phi_0+T_{U_0}(t)\Phi_0\\
	&\ +\eta_1^k\squared{-\Delta_{t,U_0}\Phi_k-\eps^2\Delta_{M}\Phi_k+\Phi_k+T_{U_0}(t)\Phi_k}+\eps^2\mathcal{R}_{U_0}[\eta_1^k,\Phi_k] \\&= \Psi+\eta_1^kb_k^\alpha\sdf{V}_\alpha.
\end{align*}
where $\mathcal{R}$ is as in \eqref{definition of R}. The above equation is satisfied if we solve the system 
\begin{align}
	-\Delta_{t,U_0}\Phi_k-\eps^2\Delta_{M}\Phi_k+\Phi_k+T_{U_0}(t)\Phi_k=-\eta_1^{-1}T_{U_0}(t)\Phi_0+b^\alpha_k\sdf{V}_\alpha\quad &\text{in }\U_{k,1}\times\R^2\label{innerprob}\\
	-\Delta_{t,U_0}\Phi_0-\eps^2\Delta_{M}\Phi_0+\Phi_0=\Psi-\eps^2\mathcal{R}_{U_0}[\eta_1^k,\Phi_k]\quad &\text{in }M\times\R^2\label{outerprob}.
\end{align}
We first solve \eqref{outerprob} with the following Lemma, which will be proved in \S\ref{Lemma positive operator in MR2}.
\begin{lemma}\label{positive operator in MR2}
	The equation 
	\begin{equation*}
		-\Delta_{t,U_0}\Phi-\eps^2\Delta_{M}\Phi+\Phi=H\quad \text{on }M\times\R^2
	\end{equation*}
	admits a solution $\Phi$ satisfying 
	\begin{equation}\label{estimate lemma6}
		\|\Phi\|_{C^{2,\gamma}(M\times\R^2)}\leq C\|H\|_{C^{0,\gamma}(M\times\R^2)}.
	\end{equation}
\end{lemma}
We use Lemma \ref{positive operator in MR2} to solve \eqref{outerprob}. Precisely, consider the two separate problems
\begin{align*}
	-\Delta_{t,U_0}\Phi_0^1-\eps^2\Delta_{M}\Phi_0^1+\Phi_0^1&=\Psi,\\
	-\Delta_{t,U_0}\Phi_0^2-\eps^2\Delta_{M}\Phi_0^2+\Phi_0^2&=-\eps^2\mathcal{R}_{U_0}[\eta_1^k,\Phi_k]
\end{align*}
and let $\boldsymbol{\Phi}=(\Phi_k)_{k\in\mathbb{N}}\in\ell^\infty(C^{2,\gamma}(\R^4))$. Then, using Lemma \ref{positive operator in MR2} we find $\Phi_0^1,\Phi_0^2$ solving the above equations and hence  $\Phi_0=\Phi_0^1+\Phi_0^2$ solving \eqref{outerprob}, satisfying 
\begin{align*}
	\|\Phi_0\|_{C^{2,\gamma}(M\times\R^2)}&\leq C(\|\Psi\|_{C^{0,\gamma}(M\times\R^2)}+\eps^2\|\mathcal{R}_{U_0}[\eta_1^k,\Phi_k]\|_{C^{0,\gamma}(M\times\R^2)})\\
	&\leq C(\|\Psi\|_{C^{0,\gamma}(M\times\R^2)}+\eps^2\|\boldsymbol{\Phi}\|_{\ell^\infty(C^{2,\gamma}(\R^4))})
\end{align*}
If we plug $\Phi_0$ in \eqref{innerprob} we obtain, for $k=1,2,\dots,$
\begin{equation*}
	-\Delta_{t,U_0}\Phi_k-\eps^2\Delta_{\xi}\Phi_k-\eps^2B_{p_k}\Phi_k+\Phi_k+T_{U_0}(t)\Phi_k=-\eta_1^{-1}T_{U_0}(t)\Phi_0+b^\alpha_k\sdf{V}_\alpha\quad \text{in }B(0,2\delta)\times\R^2,
\end{equation*}
which can be extended to the whole space by using $\eta_2^k$
\begin{equation}\label{flat formulation}
	-\Delta_{t,U_0}\Phi_k-\eps^2\Delta_{\xi}\Phi_k-\eps^2\eta_2^kB_{p_k}\Phi_k+\Phi_k+T_{U_0}(t)\Phi_k=-\eta_1^{-1}\eta_2^kT_{U_0}(t)\Phi_0+b^\alpha_k\sdf{V}_\alpha\quad \text{in }\R^n.
\end{equation}
By applying the inversion operator $\mathcal{T}$ of Lemma \ref{Solvability lin flat space}, the system \eqref{flat formulation} can be formulated as 
\begin{equation}\label{Phibold eq}
	\boldsymbol{\Phi}+\mathcal{S}(\boldsymbol{\Phi})=\mathcal{W},
\end{equation}
where 
\begin{equation*}
	\mathcal{S}(\boldsymbol{\Phi})_k=\mathcal{T}(\eta_1^{-1}\eta_2^kT_{U_0}(t)\Phi_0^2(\boldsymbol{\Phi})-\eps^2\eta_2^kB_{p_k}\Phi_k)
\end{equation*}
and 
\begin{equation*}
	\mathcal{W}_k=-\mathcal{T}(\eta_1^{-1}\eta_2^kT_{U_0}(t)\Phi_0^1(\Psi)).
\end{equation*}
We claim that $\|\mathcal{S}\|<1$, so that we can find a solution to \eqref{Phibold eq}. This concludes the existence part of Proposition \ref{invert LU0} and the sought solution $\Phi=\Phi_0+\eta_1^k\Phi_k$ satisfies 
\begin{equation*}
	\|\Phi\|_{C^{2,\gamma}(M\times\R^2)}\leq C\|\Psi\|_{C^{0,\gamma}(M\times\R^2)}.
\end{equation*}
To show the claim, we check that
\begin{align*}
	\|S(\boldsymbol{\Phi})_k\|_{C^{2,\gamma}(\R^nn)}&\leq C\|\eta_1^{-1}\eta_2^kT_{U_0}(t)\Phi_0^2(\boldsymbol{\Phi})-\eps^2\eta_2^kB_{p_k}\Phi_k\|_{C^{0,\gamma}(\R^n)}\\
	&\leq C\eps^{2}\|\boldsymbol{\Phi}\|_{\ell^\infty(C^{2,\gamma}(\R^n))},
\end{align*}
hence taking the supremum over $k$ and choosing $\eps$ sufficiently small the claim holds. To show the validity of the weighted estimates let $\varrho(t)=e^{-\sigma|t|}$ and let $\tilde{\Phi}=\varrho^{-1}\Phi$, where $\Phi$ is the solution to \eqref{corrected MR2} just found. In terms of $\tilde\Phi$ it holds 
\begin{equation*}
	-\Delta_{t,U_0}\tilde\Phi-\eps^2\Delta_{M}\tilde\Phi+\varrho^{-1}\mathcal{R}_{U_0}[\varrho,\tilde\Phi]+\tilde\Phi+T_{U_0}\tilde\Phi=\tilde\Psi+\tilde{b}^\alpha\sdf{V}_\alpha,
\end{equation*}
where $\tilde{b}^\alpha=\varrho^{-1}{b}^\alpha$, $\tilde\Psi=\varrho^{-1}\Psi$ and $\mathcal{R}$ is as in \eqref{definition of R}. Observe that 
\begin{equation*}
	\|\varrho^{-1}\mathcal{R}_{U_0}[\varrho,\tilde\Phi]\|_{C^{0,\gamma}(M\times\R^2)}\leq C\sigma\|\tilde\Phi\|_{C^{2,\gamma}(M\times\R^2)}
\end{equation*}
and thus, up to reducing $\sigma$, the invertibility theory developed in Proposition \ref{invert LU0} applies together with a fixed point argument providing the corresponding estimates for $\tilde\Phi$, namely 
\begin{equation*}
		\|\tilde\Phi\|_{C^{2,\gamma}(M\times\R^2)}\leq C\|\tilde\Psi\|_{C^{0,\gamma}(M\times\R^2)}.
\end{equation*}
The proof is concluded by observing that 
\begin{equation*}
	\|\varrho^{-1}D^2\Phi\|_{C^{0,\gamma}(M\times\R^2)}+\|\varrho^{-1}D\Phi\|_{C^{0,\gamma}(M\times\R^2)}+\|\varrho^{-1}\Phi\|_{C^{0,\gamma}(M\times\R^2)}\leq C \|\tilde\Phi\|_{C^{2,\gamma}(M\times\R^2)}.
\end{equation*}
\qed
\section{Proofs of Lemmas \ref{G contracts}, \ref{Solvability lin flat space} and \ref{positive operator in MR2}}\label{proof of lemmas}
\subsection{Proof of Lemma \ref{G contracts}}\label{G contract section}
We will prove the statement for each term $G^k$, for $k=1,2,3$. 
The claim for $G^1$ follows from a direct calculation using expression  \eqref{SW0}. 
Define 
\begin{equation*}
	\sdf{N}(h)=-\eps^{-2}\squared{\zeta_4\pi_{Z_{W,g}^\perp}N(\vphi(h))+\tilde{\sdf{B}}(h,\vphi(h))}
\end{equation*}
so that we can write 
\begin{equation*}
	G^2_\alpha(h)=\int_{\R^2}\sdf{N}(h)\cdot\sdf{V}_\alpha.
\end{equation*}
We claim that the following Lipschitz dependence hold
\begin{gather}
	\|N(\vphi_{1})-N(\vphi_{2})\|_{C^{0,\gamma}(N)}\leq C\eps^3\|\vphi_{1}-\vphi_{2}\|_{C^{2,\gamma}(N)},\label{N lip claim 1}\\
	\|\tilde{\sdf{B}}(h,\vphi_{1})-\tilde{\sdf{B}}(h,\vphi_{2})\|_{C^{0,\gamma}(N)}\leq C\eps^3\|\vphi_{1}-\vphi_{2}\|_{C^{2,\gamma}(N)},\label{N lip claim 1.5}\\
	\|\tilde{\sdf{B}}(h_{1},\vphi)-\tilde{\sdf{B}}(h_{2},\vphi)\|_{C^{0,\gamma}(N)}\leq C\eps^4\|h_1-h_2\|_{C^{2,\gamma}(M)}\|\vphi\|_{C^{2,\gamma}(N)}\label{N lip claim 2},\\
	\|\vphi_{1}-\vphi_{2}\|_{C^{2,\gamma}(N)}\leq C\eps^2\|h_1-h_2\|_{C^{2,\gamma}(M)}.\label{N lip claim 3}
\end{gather}
Assuming the validity of \eqref{N lip claim 1}-\eqref{N lip claim 3} and denoting $\vphi_l=\vphi(h_l),\ l=1,2$, we check that
\begin{align*}
		\norm{G_\alpha^2(h_1)-G_\alpha^2(h_2)}_{C^{0,\gamma}(M)}&= \norm{\int_{\R^2}[\sdf{N}(h_1)-\sdf{N}(h_2)]\cdot\sdf{V}_\alpha}_{C^{0,\gamma}(M)}\\
	&\leq C\varepsilon^{-2}\|N(\vphi_{1})-N(\vphi_{2})\|_{C^{0,\gamma}(N)}+C\varepsilon^{-2}\|B(h_1,\vphi_{1})-B(h_2\vphi_{2})\|_{C^{0,\gamma}(N)}\\
	&\leq C\eps\|\vphi_{1}-\vphi_{2}\|_{C^{2,\gamma}(N)}+ C\varepsilon^{-2}\|B(h_{1},\vphi_{1})-B(h_{1},\vphi_{2})\|_{C^{0,\gamma}(N)}\\
	&\ + C\varepsilon^{-2}\|B(h_{1},\vphi_{2})-B(h_{2},\vphi_{2})\|_{C^{0,\gamma}(N)}\\
	&\leq C\eps\|\vphi_{1}-\vphi_{2}\|_{C^{2,\gamma}(N)}+C\eps^2\|h_1-h_2\|_{C^{2,\gamma}(M)}\|\vphi\|_{C^{2,\gamma}(N)}\\
	&\leq C(\eps^3+\eps^5)\|h_1-h_2\|_{C^{2,\gamma}(M)},
\end{align*}
which proves the Lipschitz dependence of $G^2$. The validity of \eqref{N lip claim 1}, follows directly from the definition of $N$ as a second order expansion of $S$ in the direction of $\vphi$. Precisely, it holds
\begin{align*}
	\|N(\vphi_1)-N(\vphi_2)\|_{C^{0,\gamma}(N)}&\leq C \round{\|\vphi_1\|_{\infty}+\|\vphi_2\|_{\infty}}\|\vphi_1-\vphi_2\|_{C^{2,\gamma}(N)}\\
	&\leq C\eps^3\|\vphi_1-\vphi_2\|_{C^{2,\gamma}(N)}
\end{align*} 
where we used the fact that $\|\vphi_\ell\|_{C^{2,\gamma}(N)}=O(\eps^3)$.
Next we prove \eqref{N lip claim 3} and to do so we recall that $\vphi$ solves \eqref{corrected}. The term $\sdf{R}(h)=-S(W)$ has a Lipschitz dependence 
\begin{equation*}
	\|\sdf{R}(h_1)-\sdf{R}(h_2)\|_{C^{0,\gamma}(N)}\leq C\eps^2\|h_1-h_2\|_{C^{2,\gamma}(M)}
\end{equation*}
as shown by calculations analogous to those for the Lipschitz behaviour of $G^1$. The estimates from Proposition \ref{inver L_W} yield to 
\begin{align*}
	\|\vphi_1-\vphi_2\|_{C^{2,\gamma}(N)}&\leq C\|\sdf{R}(h_1)-\sdf{R}(h_2)\|_{C^{0,\gamma}(N)}+\|N(\vphi_1)-N(\vphi_2)\|_{C^{0,\gamma}(N)}\\
	&\leq C\eps^2\|h_1-h_2\|_{C^{2,\gamma}(M)}+C\eps^3\|\vphi_1-\vphi_2\|_{C^{2,\gamma}(N)}
\end{align*}
and hence, up to choosing $\eps$ sufficiently small, we have \eqref{N lip claim 3}. Then, we observe that \eqref{N lip claim 1.5} Is a direct consequence of the smallness in $\eps$ of all the coefficients of the second order operator $B$, while \eqref{N lip claim 2} is a consequence of the mild Lipschitz dependence in $h$ of such coefficients.
Finally, we prove that the term
\begin{equation*}
	G^3_\alpha(h_1)=-{\eps^{-2}}\int_{\R^2}\zeta_4L_{U_0}[\vphi]\cdot\sdf{V}_\alpha
\end{equation*}
is Lipschitz with a constant that is exponentially small in $\eps$. This is straightforward using the self-adjointness of $L_{U_0}$ and the fact that $\sdf{V}_\alpha\in\ker L_{U_0}$. Indeed, integrating by parts yields to 
\begin{align*}
	\int_{\mathbb{R}^{2}}\zeta_{4}L_{U_0}\left[\vphi\right]\cdot\mathsf{V}_{\alpha}dt&=\int_{\mathbb{R}^{2}}\zeta_{4}\left(-\Delta_{t,U_0}\vphi-\Delta_{M_{\varepsilon}}\vphi+\vphi+T_{U_0}\vphi\right)\cdot\mathsf{V}_{\alpha}dt\\&=-\int_{\mathbb{R}^{2}}(\Delta_{t}\zeta_{4})\mathsf{V}_{\alpha}\cdot\vphi dt-2\int_{\mathbb{R}^{2}}\left(\nabla_{t}\zeta_{4}\cdot\nabla_{t,U_0}\mathsf{V}_{\alpha}\right)\vphi dt
\end{align*}
where we also used the orthogonality between $\vphi$ and $\sdf{V}_\alpha$. Using the exponential decay of $\sdf{V}_\alpha$ and the usual argument with the support of $\zeta_4$ we obtain 
\begin{align*}
	\|G_\alpha^3(h_1)-G_\alpha^3(h_2)\|_{C^{0,\gamma}(N)}&\leq C\eps^{-2}e^{-\frac\delta\eps}\|\vphi_1-\vphi_2\|_\infty\\
	&\leq Ce^{-\frac\delta\eps}\|h_1-h_2\|_{C^{2,\gamma}(M)}.
\end{align*}
Finally, we can put all the estimates just found together to obtain that $G$ has an $O(\eps)$-Lipschitz constant, concluding the proof.
\qed
\subsection{Proof of Lemma \ref{Solvability lin flat space}}
We claim that for every $\Psi\in L^2(\R^n)$ we can find a solution to the system 
\begin{equation}\label{scaled version}
\begin{cases}
	-\Delta_{t,U_0}\Phi-\eps^2\Delta_{y}\Phi+\Phi+T_{U_0}\Phi=\Psi+b^\alpha(y)\sdf{V}_\alpha(t)&\text{on }\R^n\\
	\int_{\R^2}\Phi(y,t)\cdot\sdf{V}_\alpha(t)dt=0&\alpha=1,2.
\end{cases}	
\end{equation}
This result has already been established in \cite[Proposition 4]{Badran-delPino2022}. We sketch the proof for completeness. 
The differential equation in \eqref{scaled version} reads explicitly, denoting $\Phi=(\phi,\,\omega_tdt+\omega_ydy)^T$ and $\Psi=(\psi,\,\eta_tdt+\eta_ydy)^T$,
\begin{equation}\label{system broken}
\begin{cases}
-\Delta_{t}^{{A_0}}\phi-\eps^2\Delta_{y}\phi-\frac{1}{2}(1-3\left|{u_0}\right|^{2})\phi+2i\nabla^{{A_0}}{u_0}\cdot\omega_{t}=\psi+b^\alpha\sdf{V}_\alpha^1\\
-\Delta_{t}\omega_t-\eps^2\Delta_{y}\omega_t+\left|{u_0}\right|^{2}\omega_{t}+2\bracket{\nabla^{{A_0}}{u_0},i\phi}=\eta_{t}+b^\alpha\sdf{V}_\alpha^2\\
-\Delta_{t}\omega_y-\eps^2\Delta_{y}\omega_y+\left|{u_0}\right|^{2}\omega_{y}=\eta_{y}
\end{cases}
\end{equation}
where we have that the last equation is not coupled with the first two. Therefore, we need to solve two separate problems:
\begin{equation}\label{flat system}
	\begin{cases}
		-\Delta_{t}^{{A_0}}\phi-\eps^2\Delta_{y}\phi-\frac{1}{2}(1-3\left|{u_0}\right|^{2})\phi+2i\nabla^{{A_0}}{u_0}\cdot\omega_t=\psi+b^\alpha\sdf{V}_\alpha^1\\
		-\Delta_{t}\omega_t-\eps^2\Delta_{y}\omega_t+\left|{u_0}\right|^{2}\omega+2\bracket{\nabla^{{A_0}}{u_0},i\phi}=b^\alpha\sdf{V}_\alpha^2
	\end{cases}
\end{equation}
where $\omega=\omega(x,y)dx$, and the second one is 
\begin{equation}\label{uncoupled equation}
	-\Delta_{t}\omega_y-\eps^2\Delta_{y}\omega_y+f^2\omega_y=\eta_y\quad\text{in }\R^n.
\end{equation}
On \eqref{flat system} we apply Fourier transform on the $y$ variable; the transformed system reads 
\begin{equation}\label{compact-system}
	\left(\sdf{L}+\eps^2\left|\xi\right|^{2}\right)\hat\Phi=\hat\Psi+\hat b^\alpha\sdf{V}_\alpha\quad\text{in }\mathbb{R}^{n}
\end{equation}
where $\sdf{L}$ is the two-dimensional linearized around ${U_0}(t)$ and $\xi$ is the Fourier variable. Using the coercivity of $\sdf{L}$ given by \eqref{coercivity} we solve system \eqref{compact-system}, for any $\xi\in\R^{n-2}$ fixed, by Lax-Milgram theorem. Remark that this is the point in which the correction $\hat b^\alpha\sdf{V}_\alpha$ is needed since the coercivity is only true under the orthogonality condition with $\sdf{V}_\beta$, $\beta=1,2$. Applying the inverse Fourier transform we then obtain a solution to \eqref{flat system}. Equation \eqref{uncoupled equation} is solved in a similar, easier manner. This procedure also yields to 
\begin{equation*}
	\left\Vert \Phi\right\Vert _{H^{1}_{U_0}\left(\mathbb{R}^{n}\right)}\leq C\left\Vert \Psi\right\Vert _{L^{2}\left(\mathbb{R}^{n}\right)},
\end{equation*}
from which is simply proven, using elliptic interior regularity, the validity of the H\"older estimates 
\begin{equation}\label{holder flat}
	\norm{\Phi}_{C^{2,\gamma}(\R^n)}\leq C\norm{\Psi}_{C^{0,\gamma}(\R^n)}
\end{equation}
which concludes the proof.
\qed
\subsection{Proof of Lemma \ref{positive operator in MR2}}\label{Lemma positive operator in MR2}
The proof goes along the same lines to that in \cite[Lemma 7]{Badran-delPino2022}.
First, we claim that the a priori estimate 
\begin{equation}\label{estimate Linfty2}
	\norm{\Phi}_\infty\leq C\norm{H}_\infty
\end{equation}
holds for all sufficiently small $\eps>0$. By contradiction, suppose that there is a collection $\{\eps_n\}_{n\in\mathbb{N}}$ such that $\eps_n\to 0$ and a sequence of solutions 
	\begin{equation*}
		-\Delta_{t,{U_0}}\Phi_n-\eps_n^2\Delta_{M}\Phi_n+\Phi_n=H_n\quad \text{on }M\times\R^2.
	\end{equation*}
	such that
	\begin{equation*}
		\norm{\Phi_n}_\infty=1\quad\norm{H_n}_\infty\to 0.
	\end{equation*}
	Let $(y_n,t_n)\in M\times\R^2$ be  such that
	\begin{equation*}
		\abs{\Phi_n(y_n,t_n)}\geq\frac12\quad\forall n\in\mathbb{N}.
	\end{equation*}
	The compactness of $M$ implies that, passing to a subsequence, $\{y_n\}$ converges to a point $y_0\in M$. Locally around $y_0$ we can consider local coordinates for $M$ 
	\begin{equation*}
		\xi\in B(0,1)\mapsto Y_{y_0}(\xi)\in M
	\end{equation*} 
	defined by \eqref{coordinates on M}. For $n$ large enough, we set $y_n=Y_{y_0}(\xi_n)$ and define
	\begin{align*}
		\tilde{\Phi}_n(\zeta,\tau)&=\Phi_n(\xi_n+\eps_n\zeta,t_n+\tau),\\
		\tilde{H}_n(\zeta,\tau)&=H_n(\xi_n+\eps_n\zeta,t_n+\tau).
	\end{align*}
	where $\zeta\in B(0,a)\subset\R^n$, $a$ sufficiently small, and we left implicit the composition with $Y_{y_0}$. Now, using \eqref{expansion laplace-beltrami} we can see that $\tilde\Phi_n$ satisfies locally an equation of the form 
	\begin{equation*}
				-\Delta_{t,{U_0}}\tilde\Phi_n-\Delta_{\zeta}\tilde\Phi_n+o(1)D^2_\zeta\tilde\Phi_n+o(1)D_\zeta\tilde\Phi_n+\tilde\Phi_n=\tilde H_n\quad \text{in }B(0,a)\times\R^2
	\end{equation*}
	hence, passing to a subsequence we get uniform convergence over compact subsets of $B(0,a)\times\R^2$ to a bounded solution of 
	\begin{equation*}
		-\Delta_{t,{U_0}}\tilde\Phi-\Delta_{\zeta}\tilde\Phi+\tilde\Phi=0
	\end{equation*}
but the positivity of the operator implies $\tilde\Phi=0$, a contradiction.
 Now, the existence of a solution is proved using the positivity of the operator on an enlarging sequence of nested geodesic balls, along with a diagonal argument. Finally, using Schauder estimates and the a priori estimate \eqref{estimate Linfty2} we obtain \eqref{estimate lemma6}.
\qed
\appendix
\section{Construction of the outer function $\psi$}\label{A1}
We prove Lemma \ref{existence psi}, namely the existence of a smooth extension 
\begin{equation*}
	\psi : N\setminus M_h\to S^1,
\end{equation*}
to
\begin{equation*}
	\psi(x) =\frac{t}{|t|},\quad x=X_h(y,t)\in\mathrm{supp}\,\zeta_3,
\end{equation*}
where $\zeta_3$ is defined in \eqref{cutoffs}.
First, we state a useful Lemma. 
\begin{lemma}\label{tietze}
	Let $N$ be a smooth manifold and $F\subset N$ a closed submanifold with or without boundary. Let $z_0\in S^1$ and assume the existence of a smooth function 
	\begin{equation*}
		g:F\to S^1\setminus\{z_0\}.
	\end{equation*}
	Then $g$ extends smoothly to $N$, namely it exists a smooth function 
	\begin{equation*}
		\bar g:N\to S^1\setminus\{z_0\}
	\end{equation*}
	such that $\bar g\vert_F=g$.
\end{lemma}
This Lemma is just a consequence of a smooth version of Tietze's extension theorem and the fact that a punctured $S^1$ is diffeomorphic to $\R$. We divide the proof of Lemma \ref{existence psi} in different steps. First, we observe that 
\begin{equation*}
	\mathrm{supp}\,\zeta_3=\left \lbrace x\in N\ :\ d(x,M)\leq 5\delta\right\rbrace
\end{equation*}
and we explicitly define $\psi$ on the set
\begin{equation*}
	\mathcal{W}=\left\lbrace x\in N\ :\ d(x,B)\leq 5\delta\right\rbrace\setminus M_h
\end{equation*} 
where $B$ is the $n-1$ submanifold given by the hypothesis (H), for which it holds $M=\partial B$. We recall that orientability of $B$ guarantees the existence of a smooth normal vector field $n\in T^\perp B$. 
\begin{figure}
\centering
\begin{tikzpicture}		


\coordinate (O) at (1,6);
\path (O) -- +(2,.5) coordinate (O1);
\path (O1) -- +(3,.5) coordinate (O2);
\path (O2) -- +(-4,1) coordinate (O3);
\path (O) -- +(-.27,+.4) coordinate (O4);

\path (O) -- +(0,.4) coordinate (A);
\path (A) -- +(2,.5) coordinate (A1);
\path (A1) -- +(3,.5) coordinate (A2);
\path (A2) -- +(-4,1) coordinate (A3);
\path (A) -- +(-.27,+.4) coordinate (A4);

\path (O) -- +(0,.1) coordinate (Ax);
\path (Ax) -- +(2,.5) coordinate (Ax1);
\path (Ax1) -- +(3,.5) coordinate (Ax2);
\path (Ax2) -- +(-4,1) coordinate (Ax3);
\path (Ax) -- +(-.27,+.4) coordinate (Ax4);

\path (O) -- +(0,-.4) coordinate (B);
\path (B) -- +(2,.5) coordinate (B1);
\path (B1) -- +(3,.5) coordinate (B2);
\path (B2) -- +(-4,1) coordinate (B3);
\path (B) -- +(-.27,+.4) coordinate (B4);

\path (O) -- +(0,-.1) coordinate (Bx);
\path (Bx) -- +(2,.5) coordinate (Bx1);
\path (Bx1) -- +(3,.5) coordinate (Bx2);
\path (Bx2) -- +(-4,1) coordinate (Bx3);
\path (Bx) -- +(-.27,+.4) coordinate (Bx4);

\def\Ang{50}
\def\fillop{.4}
\def\fillopl{.2}
\def\dropac{.4}

\coordinate (P1) at (2.5,5);
\path (P1) -- +(\Ang:4) coordinate (P2);
\path (P2) -- +(0,2) coordinate (P3);
\path (P1) -- +(0,2) coordinate (P4);

\coordinate (I1) at (3,6.1);
\path (I1) -- +(0,.3) coordinate (I2);
\path (I1) -- +(0,.5) coordinate (I3);
\path (I1) -- +(0,.8) coordinate (I4);

\path (I1) -- +(\Ang:2.47) coordinate (J1);
\path (I2) -- +(\Ang:2.47) coordinate (J2);
\path (I3) -- +(\Ang:2.47) coordinate (J3);
\path (I4) -- +(\Ang:2.47) coordinate (J4);

\path (I4) --+ (\Ang:1.5) coordinate (IJ1h);
\path (IJ1h) --+ (0,0.2) coordinate (IJ1);
\path (IJ1) -- +(0,-.8) coordinate (IJ2);
\path (IJ1) -- +(0,-.4) coordinate (IJ3);

\path (I1) --+ (0,0.4) coordinate (I2m);
\path (I2m) --+ (\Ang-5:0.06) coordinate (I2n);
\path (J1) --+ (0,0.4) coordinate (J2m);
\path (J2m) --+ (\Ang-5:-0.05) coordinate (J2n);


\draw[green, draw opacity=0, name path = bottom_line] (P1) to (P2) to (P3);
\draw[green, draw opacity=0, name path = bottom_curve] (B1) to[out=20, in=-90] (B2) to (A2) to[out=90, in=30] (A3) to[out=210, in=90] (A4);

\draw[ultra thin, fill=orange!50,fill opacity=\fillop, draw opacity=\dropac, name intersections={of=bottom_line and bottom_curve}] (intersection-1) to (P2) to (intersection-2) to[in=-13, out=157] (J4) to (J3) to[out=180+\Ang, in=180+\Ang] (J2) to (J1) to[out=220, in=40] (IJ2) to[out=220, in=40]  (I1) to[in=184, out=19] (intersection-1);

\begin{scope}
	\clip (I2m) to[out=40, in=220] (IJ3) to[out=40, in=220] (J2m) to (J1) to[in=40, out=220] (IJ2) to[out=220, in=40] (I1) to (I2m);
\draw[ultra thin, fill=orange, fill opacity=\fillop, draw opacity=\dropac] (I1) to (I2) to[out=\Ang, in=\Ang] (I3) to (I4) to[out=40, in=220] (IJ1) to[out=40, in=220] (J4) to (J3) to[out=180+\Ang, in=180+\Ang] (J2) to (J1) to[out=220, in=40] (IJ2) to[out=220, in=40]  (I1);
\end{scope}


\draw[ultra thin, fill= white, fill opacity=0, draw opacity=\dropac] (B) to[out=-30, in=200] (B1) to[out=20, in=-90] (B2) to[out=90, in=30] (B3) to[out=210, in=90] (B4) to[out=-90, in=150] (B);
\draw[ultra thin, fill= white, fill opacity=0, draw opacity=\dropac] (Bx) to[out=-30, in=200] (Bx1) to[out=20, in=-90] (Bx2) to[out=90, in=30] (Bx3) to[out=210, in=90] (Bx4) to[out=-90, in=150] (Bx);

\foreach \x in {2,...,38} {     
\path (Bx) -- +(0,\x/200) coordinate (Ax);
\path (Ax) -- +(2,.5) coordinate (Ax1);
\path (Ax1) -- +(3,.5) coordinate (Ax2);
\path (Ax2) -- +(-4,1) coordinate (Ax3);
\path (Ax) -- +(-.27,+.4) coordinate (Ax4);
\pgfmathsetmacro\k{(1*\x+50)/1.5}
\draw[black!\k, draw opacity=\fillopl] (Ax4) to[out=-90, in=150] (Ax) to[out=-30, in=200] (Ax1) to[out=20, in=-90] (Ax2);
}

\draw[ultra thin, fill=gray!55, fill opacity=\fillopl, draw opacity=\dropac] (B4) to[out=-90, in=150] (B) to[out=-30, in=200] (B1) to[out=20, in=-90] (B2) to (Bx2) to[out=-90, in=20] (Bx1) to[in=-30, out=200] (Bx) to[in=-90, out=150] (Bx4) to (B4);

\path (Bx4) --+ (.03,.02) coordinate (Bx4co);
\draw[ultra thin, gray!50,fill=white] (Ax4) to[out=0, in=0] (Bx4);
\draw[ultra thin] (Bx4) to[out=90,in=200] (Bx4co);
\filldraw[ultra thin, black] (Bx4) circle (.001);

\path (Bx2) --+ (-.01,.01) coordinate (Bx2co);
\draw[ultra thin, gray!50,fill=white] (Ax2) to[out=180, in=180] (Bx2);
\draw[ultra thin] (Bx2) to[out=90,in=-20] (Bx2co);
\filldraw[ultra thin, black] (Bx2) circle (.001);

\draw[red, thin,fill=blue, fill opacity=\fillop, name path=manifold] (O) to[out=-30, in=200] (O1) to[out=20, in=-90] (O2) to[out=90, in=30] (O3) to[out=210, in=90] (O4) to[out=-90, in=150] (O);


\draw[ultra thin, fill= white, fill opacity=\fillopl, draw opacity=\dropac] (Ax) to[out=-30, in=200] (Ax1) to[out=20, in=-90] (Ax2) to[out=90, in=30] (Ax3) to[out=210, in=90] (Ax4) to[out=-90, in=150] (Ax);
\draw[ultra thin, fill = gray!55, fill opacity=\fillopl, draw opacity=\dropac] (A) to[out=-30, in=200] (A1) to[out=20, in=-90] (A2) to[out=90, in=30] (A3) to[out=210, in=90] (A4) to[out=-90, in=150] (A);
\draw[ultra thin, fill opacity=\fillopl, draw opacity=\dropac] (A) to[out=-30, in=200] (A1) to[out=20, in=-90] (A2) to[out=90, in=30] (A3) to[out=210, in=90] (A4) to[out=-90, in=150] (A);

\draw[ultra thin, fill=gray!55, fill opacity=\fillopl, draw opacity=\dropac] (A4) to[out=-90, in=150] (A) to[out=-30, in=200] (A1) to[out=20, in=-90] (A2) to (Ax2) to[out=-90, in=20] (Ax1) to[in=-30, out=200] (Ax) to[in=-90, out=150] (Ax4) to (A4);


\draw[ultra thin] (P1) to (P2) to (P3) to (P4) to (P1);


\begin{scope}
	\clip (I2m) to[out=40, in=220] (IJ3) to[out=40, in=220] (J2m) to (J4) to[in=40, out=220] (IJ1) to[out=220, in=40] (I4) to (I2m);
\draw[ultra thin, fill=orange, fill opacity=\fillop, draw opacity=\dropac] (I1) to (I2) to[out=\Ang, in=\Ang] (I3) to (I4) to[out=40, in=220] (IJ1) to[out=40, in=220] (J4) to (J3) to[out=180+\Ang, in=180+\Ang] (J2) to (J1) to[out=220, in=40] (IJ2) to[out=220, in=40]  (I1);
\end{scope}

\draw[blue] (I2m) to[out=40, in=220] (IJ3) to[out=40, in=220] (J2m);


\draw[ultra thin, fill=orange!50,fill opacity=\fillop, draw opacity=\dropac, name intersections={of=bottom_line and bottom_curve}] (P1) to (intersection-1) to[out=184, in=19] (I1) to (I2) to[out=\Ang, in=\Ang] (I3) to (I4) to[out=40, in=220] (IJ1) to[out=40, in=220] (J4) to[out=-13, in=157] (intersection-2) to (P3) to (P4) to (P1);


\node[blue, below] at (B3) {$B$};
\node[red, left] at (Ax4) {$M$};


\draw[->, thick] (6,9) to[out=10, in=120] (9,8);



\coordinate (X) at (8,5);
\coordinate (Y) at (13,8);

\draw[black, thin, fill=gray!10] (Y) circle (0.735cm);
\draw[black, thin, fill=gray!10] (X) circle (0.735cm);

\draw[black, line width=1.5cm] (X) to[out=20,in=180] (Y);
\draw[orange!30, line width=1.45cm] (X) to[out=20,in=180] (Y);

\draw[line width=0.025cm] (X) -- ++(-70:0.75cm);
\draw[line width=0.025cm] (X) -- ++(110:0.75cm);
\draw[line width=0.025cm] (Y) -- ++(-90:0.75cm);
\draw[line width=0.025cm] (Y) -- ++(90:0.75cm);

\filldraw[gray!10] (X) circle (0.2cm);
\filldraw[gray!10] (Y) circle (0.2cm);

\draw [thick,domain=-70:110] plot ({8+0.22*cos(\x)}, {5+0.22*sin(\x)});
\draw [thick,domain=-90:90] plot ({13-0.22*cos(\x)}, {8+0.22*sin(\x)});

\draw[blue, ultra thick] (X) to[out=20,in=180] (Y);	

\filldraw[red] (X) circle (0.05cm);
\filldraw[red] (Y) circle (0.05cm);

\node[blue] at (12,7.5) {$B$};
\node[red, right] at (Y) {$M$};
		
\node[right] at (11,6) {$\mathcal{W}$};
\node[orange] at (8,7) {$\mathcal{T}$};
\draw[thin, orange] (8.2,6.8) to (8.5,5.7);
\draw[<->] (11,7.27) -- (10.6,7.9);
\node at (11.05,7.7) {$5\delta$};

\draw[->] (9.57,6) -- +(127:0.5cm);
\node[above] at (9.57,6.1) {$n$};

\end{tikzpicture}

\caption{A cross-section of $B$ and $M=\partial B$, with a representation of the set $\mathcal{W}$ with its subset $\mathcal{T}$.}
\end{figure}
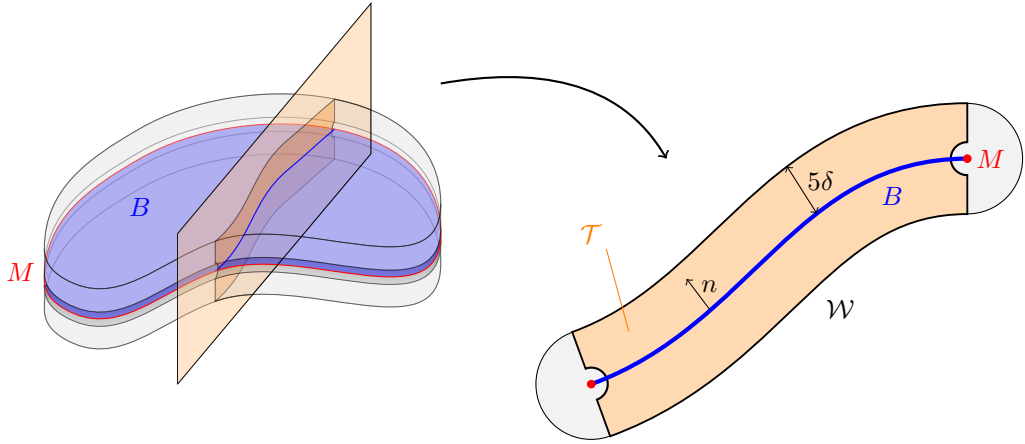
We define the sets
\begin{equation*}
	\mathcal{T}=\left\lbrace \exp_q(zn(q))\ :\ |z|\leq 5\delta,\  q\in B\right\rbrace\setminus \left\lbrace x\in \mathcal{W}\ :\ d(x,M)<\delta\right\rbrace
\end{equation*}
and
\begin{equation*}
	\mathcal{T}_\pm =\left\lbrace \exp_q(\pm5\delta n(q))\ :\ q\in B\right\rbrace.
\end{equation*}
By setting $\psi_0(x)=\pm i$ for $x\in{\mathcal{T}_\pm}$ and $\psi_0(x)=t/|t|$ in $\mathrm{supp}\zeta_3$, we get a smooth function 
\begin{equation*}
	\psi_0:(\mathcal{T}\cap\mathrm{supp}\,\zeta_3)\cup\mathcal{T}_+\cup\mathcal{T}_-\to \braces{z\in S^1\ :\ \RE{z}\geq0}
\end{equation*}
and using Lemma \ref{tietze} we find a smooth extension of $\psi_0$ on $\mathcal{T}$, namely
\begin{equation*}
	\psi_1:\mathcal{T}\to \braces{z\in S^1\ :\ \RE{z}\geq0}.
\end{equation*}
and observing that on $\mathcal{W}\setminus\mathcal{T}\subset \mathrm{supp}\zeta_3$ we can extend $\psi_1$ to the whole $\mathcal{W}$. Remark that this is already an extension of $t/|t|$. 
Next, we observe that $\psi_1^{-1}(\{1\})$ is a smooth manifold in $\mathcal{W}$ and also 
\begin{equation*}
	\mathcal{U}\coloneqq\braces{x\in N\ :\ d(x,\psi_1^{-1}(\{1\}))<\delta}\subset \mathcal{W}.
\end{equation*}
Applying again Lemma \ref{tietze} to 
\begin{equation*}
	\psi_1:\mathcal{W}\setminus\mathcal{U}\to S^1\setminus\{1\}
\end{equation*}
we can extend it to a smooth function 
\begin{equation*}
	\psi_2:N\setminus\mathcal{U}\to S^1\setminus\{1\}.
\end{equation*}
Finally, the sought function is given by
\begin{equation*}
	\psi(x)=\begin{cases}
		\psi_1(x)\quad x\in \mathcal{U}\\
		\psi_2(x)\quad x\in N\setminus\mathcal{U}.
	\end{cases}
\end{equation*}
\qed

{\bf Acknowledgements:}
	The authors have been supported  by Royal Society Research Professorship RP-R1-180114, United Kingdom.

\bibliography{Bibliography} 
\bibliographystyle{siam-fr}

\end{document}